\documentclass[11pt]{article} % perhaps use twoside option
\usepackage{geometry}
\usepackage{amsmath,amssymb,amsfonts} % need to fix margins if changing font siz
\usepackage[all]{xy}
\newcommand{\CC}{\mathbb{C}}
\newcommand{\EE}{\mathbb{E}}

\newcommand{\PP}{\mathbb{P}}

\newcommand{\ZZ}{\mathbb{Z}}
\newcommand{\Ay}{\mathcal{A}}

\newcommand{\Oh}{\mathcal{O}}
\newcommand{\De}{\mathcal{D}}
\newcommand{\eE}{\mathcal{E}}
\newcommand{\eF}{\mathcal{F}}

\newcommand{\eL}{\mathcal{L}}
\newcommand{\eM}{\mathcal{M}}
\newcommand{\Pe}{\mathcal{P}}

\newcommand{\Eff}{\mathfrak{E}}
\newcommand{\HH}{H\!H}
\newcommand{\NHH}{N\!H\!H}
\newcommand{\Perf}{\mathrm{Perf}}
\newcommand{\al}{\alpha}
\newcommand{\be}{\beta}
\newcommand{\gam}{\gamma}
\newcommand{\Del}{\Delta}

\newcommand{\fie}{\varphi}
\newcommand{\Fie}{\Phi}
\newcommand{\Abar}{\overline{A}}
\newcommand{\Bbar}{\overline{B}}
\newcommand{\Cbar}{\overline{C}}

\newcommand{\Ebar}{\overline{E}}
\newcommand{\Gam}{\Gamma}
\newcommand{\Gambar}{\overline\Gamma}

\newcommand{\La}{\Lambda}
\newcommand{\Lchi}{\eL_{\chi}}
\newcommand{\kbar}{\overline k}
\newcommand{\id}{\mathrm{id}}

\newcommand{\lin}{\sim}
\newcommand{\isom}{\cong}
\newcommand{\Gtilde}{\widetilde G}

\newcommand{\qed}{\hfill\square}

\newcommand{\Span}[1]{\left<#1\right>}

\DeclareMathOperator{\Pic}{Pic}
\DeclareMathOperator{\Tors}{Tors}
\DeclareMathOperator{\Hom}{Hom}

\DeclareMathOperator{\Ext}{Ext}
\DeclareMathOperator{\diag}{diag}

\newtheorem{thm}{Theorem}[section]
\newtheorem{prop}{Proposition}[section]
\newtheorem{cor}{Corollary}[section]
\newtheorem{lemma}{Lemma}[section]
\newtheorem{rmk}{Remark}[section]
\newtheorem{dfn}{Definition}[section]
\newenvironment{pf}{\paragraph{Proof}}{\par\medskip}
\newenvironment{pf1}{\paragraph{First Proof}}{\par\medskip}
\newenvironment{pf2}{\paragraph{Second Proof}}{\par\medskip}

\title{Enumerating exceptional collections of line bundles\\
on some surfaces of general type}
\author{Stephen Coughlan\thanks{Former:
Department of Mathematics and Statistics
Lederle Graduate Research Tower
University of Massachusetts
Amherst, MA 01003-9305
coughlan@math.umass.edu} \thanks{Current:
Institute of Algebraic Geometry,
Leibniz Universit\"at Hannover
Welfengarten 1
30167 Hannover, Germany
coughlan@math.uni-hannover.de}}
\date{}

\begin{document}
\maketitle
\begin{abstract}
We use constructions of surfaces as abelian covers to write down exceptional collections of line bundles of maximal length for every surface $X$ in certain families of surfaces of general type with $p_g=0$ and $K_X^2=3,4,5,6,8$. We also compute the algebra of derived endomorphisms for an appropriately chosen exceptional collection, and the Hochschild cohomology of the corresponding quasiphantom category. As a consequence, we see that the subcategory generated by the exceptional collection does not vary in the family of surfaces. Finally, we describe the semigroup of effective divisors on each surface, answering a question of Alexeev.
\end{abstract}

\setcounter{tocdepth}{1}
\tableofcontents
\section{Introduction}
Exceptional collections of maximal length on surfaces of general type with $p_g=0$ have been constructed for Godeaux surfaces \cite{BBS,BBKS}, primary Burniat surfaces \cite{AO}, and Beauville surfaces \cite{GS,L}. Recently, progress has also been made for some fake projective planes \cite{minifolds,Fak}. In this article, we present a method which can be applied uniformly to produce exceptional collections of line bundles on several surfaces with $p_g=0$, including Burniat surfaces with $K^2=6$ (cf.~\cite{AO}), $5,4,3$, Kulikov surfaces with $K^2=6$ and some Beauville surfaces with $K^2=8$ \cite{GS,L}. In fact we do more: we enumerate all exceptional collections of line bundles corresponding to any choice of numerical exceptional collection. We can use this enumeration process to find those exceptional collections that are particularly well-suited to studying the surface itself, and possibly its moduli space.

Both \cite{AO} and \cite{GS} hinted that it should be possible to produce exceptional collections of line bundles on a wide range of surfaces of general type with $p_g=0$. This inspired us to build the approaches of \cite{AO,GS} into the larger framework of abelian covers (see especially Section \ref{sec!prelim}), an important part of which is a new formula for the pushforward of certain line bundles on any abelian cover, generalising formulas of Pardini \cite{Pa}. We believe that this work is a step in the right direction, even though there remain many families of surfaces which require further study (see Section \ref{sec!strategy} for more details).

Let $X$ be a surface of general type with $p_g=0$, and let $Y$ be a del Pezzo surface with $K_Y^2=K_X^2$. The lattices $\Pic X/{\Tors X}$ and $\Pic Y$ are both isomorphic to $\ZZ^{1,N}$, where $N=9-K_X^2$, and moreover, the cohomology groups $H^2(X,\ZZ)$ and $H^2(Y,\ZZ)$ are completely algebraic. By exploiting this relationship between $X$ and $Y$, we can study exceptional collections of line bundles on $X$. Indeed, exceptional collections on del Pezzo surfaces are well understood after \cite{O}, \cite{KO}, and we sometimes refer to $X$ as a fake del Pezzo surface, to emphasise this analogy.

Suppose now that $X$ is a fake del Pezzo surface that is constructed as a branched Galois abelian cover $\fie\colon X\to Y$, where $Y$ is a (weak) del Pezzo surface with $K_Y^2=K_X^2$. Many fake del Pezzo surfaces can be constructed in this way \cite{BCGP}, but we require certain additional assumptions on the branch locus and Galois group (see Section \ref{sec!strategy}). These assumptions ensure that there is an appropriate choice of lattice isometry $\Pic Y\to\Pic X/{\Tors X}$. This isometry is combined with our pushforward formula to calculate the coherent cohomology of any line bundle on $X$.
\begin{thm}[Theorem \ref{thm!pushforward}]\label{thm!intro-pushforward} Let $X$ be a fake del Pezzo surface satisfying our assumptions, and let $L$ be any line bundle on $X$. We have an explicit formula for the line bundles $M_\chi$ appearing in the pushforward $\fie_*L=\bigoplus_{\chi\in G^*} M_\chi$, where $G$ is the Galois group of the cover $\fie\colon X\to Y$.
\end{thm}

Working modulo torsion, we can use the above lattice isometry to lift any exceptional collection of line bundles on $Y$ to a numerical exceptional collection on $X$. We then incorporate Theorem \ref{thm!intro-pushforward} into a systematic computer search, to find those combinations of torsion twists which correspond to an exceptional collection on $X$.

The search for exceptional collections on fake del Pezzo surfaces, leads naturally to the following question, which was asked by Alexeev \cite{Alexeev}: \begin{quote}Can we characterise effective divisors on $X$ in terms of those on $Y$?\end{quote}
For example, in \cite{Alexeev}, Alexeev gives an explicit description of the semigroup of effective divisors on the Burniat surface with $K^2=6$, and proposes similar descriptions for the other Burniat surfaces. We use our pushforward formula to prove these characterisations for the Burniat surfaces and other fake del Pezzo surfaces, cf.~Theorems \ref{thm!Kulikov-effective}, \ref{thm!Burniat-effective}. 

\begin{thm}\label{thm!intro-effective} Let $X$ be a fake del Pezzo surface satisfying our assumptions. Then the semigroup of effective divisors on $X$ is generated by the reduced pullback of irreducible components of the branch divisor, together with pullbacks of certain $(-1)$-curves on $Y$.
\end{thm}

Let $\EE$ be an exceptional collection on $X$, and suppose $H_1(X,\ZZ)$ is nontrivial. Then $\EE$ can not be full, for $K$-theoretic reasons (see Section \ref{sec!heights}). Hence we have a semiorthogonal decomposition of the bounded derived category of coherent sheaves on $X$:
\[D^b(X)={\Span{\EE,\Ay}}.\]
If $\EE$ is of maximal length, then $\Ay$ is called a quasiphantom category; that is, $K_0(\Ay)$ is torsion and the Hochschild homology $\HH_*(\Ay)$ is trivial. Even when $H_1(X,\ZZ)$ vanishes, an exceptional collection of maximal length need not be full (see \cite{BBKS}), and in this case $\Ay$ is called a phantom category, because $K_0(\Ay)$ is trivial.

On the other hand, the Hochschild cohomology does detect the quasiphantom category $\Ay$; in fact, $\HH^*(\Ay)$ measures the formal deformations of $\Ay$. We calculate $\HH^*(\Ay)$ by considering the $A_\infty$-algebra of endomorphisms of $\EE$, together with the spectral sequence developed in \cite{K}. Indeed, one of the advantages of our systematic search, is that we can find exceptional collections for which the higher multiplications in the $A_\infty$-algebra of $\EE$ are as simple as possible. Theorem \ref{thm!main} below serves as a prototype statement of our results for a good exceptional collection on a fake del Pezzo surface. More precise statements can be found for the Kulikov surface in Section \ref{sec!Kulikov-quasiphantom}.
\begin{thm}\label{thm!main} Let $\mathcal{X}\to T$ be a family of fake del Pezzo surfaces satisfying our assumptions. Then for any $t$ in $T$, there is an exceptional collection $\EE$ of line bundles on $X=\mathcal{X}_t$ which has maximal length $12-K_X^2$. Moreover, the subcategory of $D^b(X)$ generated by $\EE$ does not vary with $t$, and the Hochschild co\-homology of $X$ agrees with that of the quasi-phantom category $\Ay$ in degrees less than or equal to two.
\end{thm}

The significance of Theorem \ref{thm!main} is amplified by the reconstruction theorem of \cite{BO}: if $X$ and $X'$ are smooth, $\pm K_X$ is ample, and $D^b(X)$ and $D^b(X')$ are equivalent bounded derived categories, then $X\isom X'$. In conjunction with Theorem \ref{thm!main}, we see that if $K_X$ is ample, then $X$ can be reconstructed from the quasi-phantom category $\Ay$. The gluing between $\Ay$ and $\EE$ does not vary with $X$, because the statement about Hochschild cohomology implies that the formal deformation spaces of $X$ are isomorphic to the formal deformation spaces of $\Ay$. Currently, it is not clear whether there is any practical way to extract information about $X$ from $\Ay$, although some interesting ideas are discussed in \cite{AO}. It would be interesting to know whether this ``rigidity'' of $\EE$ is a general phenomenon, or just a coincidence for good choices of exceptional collection.

In Section \ref{sec!prelim} we review abelian covers, and prove our result on pushforwards of line bundles, which is valid for any abelian cover, and is used throughout. In Section \ref{sec!strategy}, we explain our assumptions on the fake del Pezzo surface $X$ and its Galois covering structure $\fie\colon X\to Y$, and describe our approach to enumerating exceptional collections on the surface of general type. Section \ref{sec!Kulikov} is an extended treatment of the Kulikov surface with $K^2=6$, which is an example of a fake del Pezzo surface. We give a cursory review of dg-categories and $A_\infty$-algebras in Section \ref{sec!heights}, as background to our discussion of quasi-phantom categories and the theory of heights from \cite{K}. We then show how to compute the $A_\infty$-algebra and height of an exceptional collection on the Kulikov surface. In Section \ref{sec!burniat-secondary} we prove Theorem \ref{thm!intro-effective} for the secondary nodal Burniat surface with $K^2=4$.
Appendix \ref{app!Kulikov} lists certain data relevant to the Kulikov surface example of Section \ref{sec!Kulikov}, and Appendix \ref{app!Burniat-surface-data} applies similarly to the secondary nodal Burniat surface of Section \ref{sec!burniat-secondary}.

With appropriate amendments, Theorems \ref{thm!intro-effective} and \ref{thm!main} hold for the Burniat surfaces with $K^2=6,5,4,3$ and some Beauville surfaces with $K^2=8$. The arguments used are similar to those appearing in Sections \ref{sec!Kulikov} and \ref{sec!Burniat}, and we refer to \cite{computer-code} for details. We have exceptional collections of maximal length on the tertiary Burniat surface with $K^2=3$. In this case it is necessary to use the Weyl group action on the Picard group to find exceptional collections. We can show that the $A_\infty$-category is formal, but we do not yet know how to compute the Hochschild cohomology of the quasiphantom category.

In order to use results on deformations of each fake del Pezzo surface, we work over $\CC$.

\begin{rmk}\rm
The calculation of $\fie_*L$ according to Theorem \ref{thm!intro-pushforward} is elementary but repetitive; we include a few sample calculations to illustrate how to do it by hand, but when the torsion group becomes large, it is more practical to use computer algebra. Our enumerations of exceptional collections are obtained by simple exhaustive computer searches. We use Magma \cite{Magma}, and the annotated scripts are available from \cite{computer-code}.
\end{rmk}
\paragraph{Acknowledgements} I would like to thank Valery Alexeev, Ingrid Bauer, Gavin Brown, Fabrizio Catanese, Paul Hacking, Al Kasprzyk, Anna Kazanova, Alexander Kuznetsov, Miles Reid and Jenia Tevelev for helpful conversations or comments about this work. I thank the DFG for support during part of this work through grant Hu 337-6/2.
\section{Preliminaries}\label{sec!prelim}
We collect together the relevant material on abelian covers. See especially \cite{Pa}, \cite{BC} or \cite{Ku} for details. Unless stated otherwise, $X$ and $Y$ are normal projective varieties, with $Y$ nonsingular. Let $G$ be a finite abelian group acting faithfully on $X$ with quotient $\fie\colon X\to Y$. Write $\Del=\sum\Del_i$ for the branch locus of $\fie$, where each $\Del_i$ is a reduced, irreducible effective divisor on $Y$. The cover $\fie$ is determined by the group homomorphism
\[\Fie\colon\pi_1(Y-\Del)\to H_1(Y-\Del,\ZZ)\to G,\]
which assigns an element of $G$ to the class of a loop around each irreducible component $\Del_i$ of $\Del$. If $\Fie$ is surjective, then $X$ is irreducible. The factorisation through $H_1(Y-\Del,\ZZ)$ arises because $G$ is assumed to be abelian, so we only need to consider the map $\Fie\colon H_1(Y-\Del,\ZZ)\to G$. For brevity, we refer to the loop around $\Del_i$ by the same symbol, $\Del_i$.

Let $\widetilde Y$ be the blow up of $Y$ at a point $P$ where several branch components $\Del_{i_1},\dots,\Del_{i_k}$ intersect. Then there is an induced cover of $\widetilde Y$, and the image of the exceptional curve $E$ under $\Fie$ is given by
\begin{equation}\label{eqn!exceptional-curve}\Fie(E)=\sum_{j=1}^k\Fie(\Del_{i_j}).\end{equation}

Fix an irreducible reduced component $\Gam$ of $\Del$ and denote $\Fie(\Gam)$ by $\gam$. Then the inertia group of $\Gam$ is the cyclic group $H\subset G$ generated by $\gam$. Choosing the generator of $H^*=\Hom(H,\CC^*)$ to be the dual character $\gam^*$, we may identify $H^*$ with $\ZZ/n$, where $n$ is the order of $\gam$. Composing the restriction map $\text{res}\colon G^*\to H^*$ with this identification gives
\[G^*\to\ZZ/n,\ \chi\mapsto k,\]
where $\chi|_H=(\gam^*)^k$ for some $0\le k\le n-1$. On the other hand, given $\chi$ in $G^*$ of order $d$, the evaluation map $\chi\colon G\to \ZZ/d$ satisfies
\[\chi(\gam)=\tfrac dn\chi|_H(\gam)=\tfrac {dk}n\]
as a residue class in $\ZZ/d$ (or as an integer between $0$ and $d-1$).

The pushforward of $\fie_*\Oh_X$ breaks into a direct sum of eigensheaves
\begin{equation}\label{eqn!charsheaves}
\fie_*\Oh_X=\bigoplus_{\chi\in G^*}\Lchi^{-1}.
\end{equation}
Moreover, the $\Lchi$ are line bundles on $Y$ and by Pardini \cite{Pa}, their associated (integral) divisors $L_\chi$ are given by the formula
\begin{equation}\label{eqn!pardini}dL_\chi=\sum_i\chi\circ\Fie(\Del_i)\Del_i.\end{equation}
The line bundles $\Lchi$ play a pivotal role in the sequel, and we refer to them as the \emph{character sheaves} of the cover $\fie\colon X\to Y$.

\subsection{Line bundles on $X$}

We develop tools for calculating with torsion line bundles on $X$. Let $\pi'\colon A'\to X$ be the maximal abelian cover of $X$; that is, the \'etale cover of $X$ associated to the subgroup $\pi_1(X)^{\text{ab}}=H_1(X,\ZZ)$ of $\pi_1(X)$. Now let $\psi'$ be the composite map $\fie\circ\pi'\colon A'\to Y$. It is not always true that $\psi'$ is Galois and ramified over the same branch divisor $\Del$ as $\fie\colon X\to Y$ (see for example \cite{Shab}, \cite{BCF}). So choose a maximal subgroup $T$ of the torsion subgroup $\Tors X$ in $\Pic X$ whose associated cover $\psi\colon A\to Y$ is Galois and ramified over $\Del$. We have the following commutative diagram
\[\xymatrix%@R=30pt@C=35pt
{&A\ar[dl]_{\pi}\ar[dr]^{\psi}&\\
X\ar[rr]^{\fie}&&Y}\]
Let the Galois group of $\psi$ be $\Gtilde$. Then the original group $G$ is the quotient $\Gtilde/T$, so we get short exact sequences
\begin{equation}\label{eqn!sesgroup}0\to T\to\Gtilde\to G\to 0\end{equation}
and
\begin{equation}\label{eqn!sesdual}0\leftarrow T^*\leftarrow \Gtilde^*\leftarrow G^*\leftarrow0\end{equation}
where $G^*=\Hom(G,\CC^*)$, etc. In fact, for each surface that we consider, these exact sequences are split, so that
\begin{equation}\label{eqn!directsum}\Gtilde=G\oplus T,\ \Gtilde^*=G^*\oplus T^*.\end{equation}

Let $\Gambar$ be a reduced irreducible component of the branch locus $\Del$ of an abelian cover $\fie\colon X\to Y$ and suppose the inertia group of $\Gambar$ is cyclic of order $n$. Then
\begin{dfn}[cf.~\cite{AO}]\label{def!redpull} The \emph{reduced pullback} $\Gam$ of $\Gambar$ is the (integral) divisor $\Gam=\frac1n\fie^*(\Gambar)$ on $X$.
\end{dfn}

\begin{rmk}\rm The reduced pullback extends to arbitrary linear combinations $\sum_ik_i\Del_i$ in the obvious way. We use a bar to denote divisors on $Y$ and remove the bar when taking the reduced pullback. In other situations, it is convenient to use $D_i$ to denote the reduced pullback of a branch divisor $\Del_i$.
\end{rmk}

The remainder of this section is dedicated to calculating the pushforward $\fie_*(L\otimes\tau)$, where $L=\Oh_X(\sum_i k_iD_i)$ is the line bundle associated to the reduced pullback of $\sum_i k_i\Del_i$, and $\tau$ is any torsion line bundle contained in $T\subset\Tors X$. We do this by exploiting the association of the free part $L$ with $\fie\colon X\to Y$, and the torsion part $\tau$ with $\pi\colon A\to X$. The formulae that we obtain are a natural extension of results in \cite{Pa}. It may be helpful to skip ahead to Examples \ref{exa!Campedelli} and \ref{exa!Campedelli-cont} before reading this section in detail.

\subsection{Free case}\label{sec!free}

Until further notice, we write $\Gambar\subset Y$ for an irreducible component of the branch divisor $\Del$ of $\fie\colon X\to Y$. By Pardini \cite{Pa}, the inertia group $H\subset G$ of $\Gambar$ is cyclic, and $H$ is generated by $\Fie(\Gambar)$ of order $n$. Let $\Gam\subset X$ be the reduced pullback of $\Gambar$, so that $n\Gam=\fie^*(\Gambar)$. We start with cyclic covers.

\begin{lemma}\label{lem!cyclic} Let $\al\colon X\to Y$ be a cyclic cover with group $H\isom\ZZ/n$, and suppose that $\Gambar$ is an irreducible reduced component of the branch divisor. Let $\Gam$ be the reduced pullback of $\Gambar$, and suppose $0\le k\le{n-1}$. Then
\[\al_*\Oh_X(k\Gam)=\bigoplus_{i\in H^*-S}\eM_i^{-1}\oplus\bigoplus_{i\in S}\eM_i^{-1}(\Gambar),\]
where $\eM_i$ is the character sheaf associated to $\al$ with character $i\in H^*$, and \[S=\{n-k,\dots,n-1\}\subset H^*\isom\ZZ/n.\]
\end{lemma}

\begin{rmk}\label{rmk!projformula}\rm If $k$ is a multiple of $n$, say $k=pn$, the projection formula gives
\[\al_*\Oh_X(k\Gam)=\al_*(\al^*\Oh_Y(p\Gambar))=\al_*\Oh_X\otimes\Oh_Y(p\Gambar)=\bigoplus_{i\in H^*}\eM_i^{-1}(p\Gambar).\]
Thus the lemma extends to any integer multiple of $\Gam$.
\end{rmk}

\begin{pf}
After removing a finite number of points from $\Gambar$, we may choose a neighbourhood $U$ of $\Gambar$ such that $U$ does not intersect any other irreducible components of $\Del$. Then since $X$ and $Y$ are normal we may calculate $\al_*\Oh_X(k\Gam)$ locally on $\al^{-1}(U)$ and $U$. In what follows, we do not distinguish $U$ (respectively $\al^{-1}(U)$) from $Y$ (resp.~$X$).

Let $g=\Fie(\Gambar)$ so that $H=\Span{g}\isom\ZZ/n$, and identify $H^*$ with $\ZZ/n$ via $g^*=1$. Locally, write $\al\colon\al^{-1}(U)\to U$ as $z^n=b$ where $b=0$ defines $\Gambar$ in $U$. Then
\[\al_*\Oh_X=\bigoplus_{i=0}^{n-1}\Oh_Yz^i=\bigoplus_{i=0}^{n-1}\Oh_Y(-\tfrac in\Gambar)=\bigoplus_{i=0}^{n-1}\eM_i^{-1},\]
where the last equality is given by \eqref{eqn!pardini}. Thus $\al_*\Oh_X$ is generated by $1,z,\dots,z^{n-1}$ as an $\Oh_Y$-module, and the $\Oh_Y$-algebra structure on $\al_*\Oh_X$ is induced by the equation $z^n=b$.

The calculation for $\Oh_X(k\Gam)$ is similar,
\[\al_*\Oh_X(k\Gam)=\al_*\Oh_X\frac1{z^k}=\bigoplus_{i=-k}^{n-k-1}\Oh_Yz^i=\bigoplus_{i=0}^{n-k-1}\Oh_Yz^i\oplus\bigoplus_{i=-k}^{-1}\Oh_Y\frac {z^{n+i}}b\]
where we use $z^n=b$ to remove negative powers of $z$. Thus
\begin{align*}\al_*\Oh_X(k\Gam)&=\bigoplus_{i=0}^{n-k-1}\Oh_Y(-\tfrac in\Gambar)\oplus\bigoplus_{i=n-k}^{n-1}\Oh_Y(-\tfrac in\Gambar)(\Gambar)\\
&=\bigoplus_{i\in H^*-S}\eM_i^{-1}\oplus\bigoplus_{i\in S}\eM_i^{-1}(\Gambar),
\end{align*}
where $S=\{n-k,\dots,n-1\}$.
$\qed$
\end{pf}

The lemma can be extended to any abelian group using arguments inspired by Pardini \cite{Pa} Sections 2 and 4.

\begin{prop}\label{prop!pushfree} Let $\fie\colon X\to Y$ be an abelian cover with group $G$, and let $k=np+\kbar$, where $0\le\kbar\le n-1$. Then
\[\fie_*\Oh_X(k\Gam)=\bigoplus_{\chi\in G^*-S_{k\Gambar}}\Lchi^{-1}(p\Gambar)\oplus\bigoplus_{\chi\in S_{k\Gambar}}\Lchi^{-1}((p+1)\Gambar),\]
where
\[S_{k\Gambar}=\{\chi\in G^*: n-\kbar\le\chi|_{H}\le n-1\}.\]
\end{prop}

\begin{pf}
By the projection formula, we only need to consider the case $k=\kbar$ (cf.~Remark \ref{rmk!projformula}). As in the proof of Lemma \ref{lem!cyclic}, after removing a finite number of points, we may take a neighbourhood $U$ of $\Gambar$ which does not intersect any other components of $\Del$. We work on $U$ and its preimages $\fie^{-1}(U)$, $\be^{-1}(U)$.

Factor $\fie\colon X\to Y$ as
\[X\xrightarrow{\al}Z\xrightarrow{\be}Y,\]
where $\al$ is a cyclic cover ramified over $\Gam$ with group $H=\Span{g}\isom\ZZ/n$, and $\be$ is unramified by our assumptions. As in Lemma \ref{lem!cyclic} we denote the character sheaves of $\al$ by $\eM_i$, and those of the composite map $\fie=\be\circ\al$ by $\Lchi$. Now
\begin{equation}\label{eqn!unram}\be_*\eM_i=\bigoplus_{\chi\in[i]}\Lchi\end{equation}
where the notation $[i]$ means the preimage of $i$ in $H^*$ under the restriction map $\text{res}\colon G^*\to H^*$. That is,
\[[i]=\{\chi\in G^*:\chi|_H=i\},\]
where $d$ is the order of $\chi$. Since $\be$ is not ramified we combine Lemma \ref{lem!cyclic} and \eqref{eqn!unram} to get
\[\fie_*\Oh_X(k\Gam)=\bigoplus_{\chi\in G^*-S_{k\Gambar}}\Lchi^{-1}\oplus\bigoplus_{\chi\in S_{k\Gambar}}\Lchi^{-1}(\Gambar)\]
where 
\[S_{k\Gambar}=\{\chi\in G^* : n-k\le\chi|_H\le n-1\}\]
is the preimage of $S=\{n-k,\dots,n-1\}\subset H^*$ under $\text{res}\colon G^*\to H^*$.
$\qed$
\end{pf}
\subsubsection{Example (Campedelli surface)}\label{exa!Campedelli}
Let $\fie\colon X\to\PP^2$ be a $G=(\ZZ/2)^3$-cover branched over seven lines in general position. We label the lines $\Del_1,\dots,\Del_7$, and define $\Fie$ to induce a set-theoretic bijection between $\{\Del_i\}$ and $(\ZZ/2)^3-\{0\}$. We make the definition of $\Fie$ more precise later (see Example \ref{exa!Campedelli-cont}). It is well known (\cite[\S 4]{Ku}) that $X$ is a surface of general type with $p_g=0$, $K^2=2$ and $\pi_1=(\ZZ/2)^3$.

Choose generators $g_1,g_2,g_3$ for $(\ZZ/2)^3$ so that $\Fie(\Del_1)=g_1$. There are eight character sheaves for the cover, which we calculate using formula \eqref{eqn!pardini},
\[\eL_{(0,0,0)}=\Oh_{\PP^2},\ \eL_{\chi}=\Oh_{\PP^2}(2)\text{ for }\chi\ne(0,0,0).\]
Write $D_1$ for the reduced pullback of $\Del_1$, so that $\fie^*(\Del_1)=2D_1$. Then 
\[S_{\Del_1}=\{\chi : \chi|_{\Span{g_1}}=1\}=\{(1,0,0),(1,1,0),(1,0,1),(1,1,1)\},\]
so that by Proposition \ref{prop!pushfree}, we have
\[\fie_*\Oh_X(D_1)=\Oh_{\PP^2}\oplus4\Oh_{\PP^2}(-1)\oplus3\Oh_{\PP^2}(-2).\]

\subsection{Torsion case}\label{sec!tors}

In this section we use the maximal abelian cover $A$ to calculate the pushforward of a torsion line bundle on $X$. To simplify notation, we assume that the composite cover $A\to X\to Y$ is Galois with group $\Gtilde$, so that $T=\Tors X$.

\begin{prop}\label{prop!torsion} Let $\tau$ be a torsion line bundle on $X$. Then
\[\fie_*\Oh_X(-\tau)=\bigoplus_{\chi\in G^*}\eL_{\chi+\tau}^{-1}.\]
where addition $\chi+\tau$ takes place in $\Gtilde^*=G^*\oplus T^*$.
\end{prop}
\begin{rmk}\rm
Note that $\eL_{\chi+\tau}$ is a character sheaf for the $\Gtilde$-cover $\fie\colon A\to Y$, and the proposition allows us to interpret $\eL_{\chi+\tau}$ as a character sheaf for the $G$-cover $\fie\colon X\to Y$. Unfortunately, there is still some ambiguity, because we do not determine which character in $G^*$ is associated to each $\eL_{\chi+\tau}$ under the splitting of exact sequence \eqref{eqn!sesdual}. On the other hand, the special case $\tau=0$ gives
\[\fie_*\Oh_X=\bigoplus_{\chi\in G^*}\Lchi^{-1}.\]
\end{rmk}
\begin{pf} The structure sheaf $\Oh_A$ decomposes into a direct sum of the torsion line bundles when pushed forward to $X$
\[\pi_*\Oh_A=\bigoplus_{\tau\in\Tors X}\Oh_X(-\tau).\]
Thus $\Oh_X(\tau)$ is the character sheaf with character $\tau$ under the identification $T^*\isom\Tors X$. The composite $\fie_*\pi_*\Oh_A$ breaks into character sheaves according to \eqref{eqn!charsheaves}, and the image of $\Oh_X(-\tau)$ is the direct sum of those character sheaves with character contained in the coset $G^*+\tau$ of $\tau$ in $\Gtilde^*$ under \eqref{eqn!directsum}.
$\qed$
\end{pf}

\subsection{General case}\label{sec!general}

Now we combine Propositions \ref{prop!pushfree} and \ref{prop!torsion} to give our formula for pushforward of line bundles $\Oh_X(\sum_iD_i)\otimes\tau$. The formula looks complicated, but most of the difficulty is in the notation.
\begin{dfn}\label{def!S_i}\rm Let $n_i$ be the order of $\Psi(\Del_i)$ in $\Gtilde$, and write $k_i=n_ip_i+\kbar_i$, where $0\le\kbar_i\le n_i-1$. Then given any subset $I\subset\{1,\dots,m\}$, we define
\[S_I[\tau]=\bigcap_{i\in I} S_{k_i\Del_i}[\tau]\cap\bigcap_{j\in I^c} S_{k_j\Del_j}[\tau]^c,\]
where 
\[S_{k\Gambar}[\tau]=\{\chi\in G^* : n-\kbar\le\tfrac nd(\chi+\tau)(\Psi(\Gambar))\le n-1\}\]
for any reduced irreducible component $\Gambar$ of the branch locus $\Del$. Note that for fixed $\tau$ in $T^*$, the collection of all $S_I[\tau]$ partitions $G^*$.
\end{dfn}

\begin{thm}\label{thm!pushforward} Let $D=\sum_{i=1}^mk_iD_i$ be the reduced pullback of the linear combination of branch divisors $\sum_{i=1}^mk_i\Del_i$ on $Y$. Then
\[\fie_*\Oh_X(D-\tau)=\bigoplus_{I}\bigoplus_{\chi\in S_I[\tau]}\eL_{\chi+\tau}^{-1}(\Del_I),\]
where $I$ is any subset of $\{1,\dots,m\}$ and $\Del_I=\sum_{i\in I}\Del_i$.
\end{thm}
\begin{rmk}\rm For simplicity, we have assumed that $k_i=\kbar_i$ for all $i$ in the statement and proof of the theorem. When this is not the case, by the projection formula (cf.~Remark \ref{rmk!projformula}) we twist by $\Oh_Y(\sum_{i=1}^mp_i\Del_i)$.
\end{rmk}
\begin{pf}
Fix $i$ and let $D_i$ be the reduced pullback of an irreducible component $\Del_i$ of the branch divisor. Choose a neighbourhood of $\Del_i$ which does not intersect any other $\Del_j$. This may also require us to remove a finite number of points from $D_i$. We work locally in this neighbourhood and its preimages under $\fie$, $\pi$.

Now by the projection formula,
\[\pi_*\pi^*\Oh_X(k_iD_i)=\pi_*\Oh_A\otimes\Oh_X(k_iD_i),\]
and thus
\[\psi_*\pi^*\Oh_X(k_iD_i)=\bigoplus_{\tau\in T}\fie_*\Oh_X(k_iD_i-\tau).\]
Then we combine Propositions \ref{prop!pushfree} and \ref{prop!torsion} to obtain
\[\fie_*\Oh_X(k_iD_i-\tau)=\bigoplus_{\chi\in G^*-S_{k_i\Del_i}[\tau]}\eL^{-1}_{\chi+\tau}\oplus\bigoplus_{\chi\in S_{k_i\Del_i}[\tau]}\eL^{-1}_{\chi+\tau}(\Del_i),\]
where the indexing is explained in Definition \ref{def!S_i}.

To extend to the global setting and linear combinations $\sum k_iD_i$, we just need to keep track of which components of $\Del$ should appear as a twist of each $\eL^{-1}_{\chi+\tau}$ in the direct sum. This book-keeping is precisely the purpose of Definition \ref{def!S_i}.
$\qed$
\end{pf}

Using the formula
\begin{equation}\label{eqn!canonical}
K_X=\fie^*\big(K_Y+\sum_i\tfrac{n_i-1}{n_i}\Del_i\big)
\end{equation}
and the Theorem, we give an alternative proof of the decomposition of $\fie_*\Oh_X(K_X)$.
\begin{cor}{\cite[Proposition 4.1]{Pa}}\label{cor!canonical} We have
\[\fie_*\Oh_X(K_X)=\bigoplus_{\chi\in G^*}\eL_{\chi^{-1}}(K_Y).\]
\end{cor}
\begin{pf}
Let $D_i$ be the reduced pullback of $\Del_i$. Then by \eqref{eqn!canonical} and the projection formula, we have
\begin{align*}
\fie_*(\Oh_X(K_X))&=\fie_*\Big(\fie^*\Oh_Y(K_Y)\otimes\Oh_X\Big(\sum_i(n_i-1)D_i\Big)\Big)\\
&=\Oh_Y(K_Y)\otimes\fie_*\Oh_X\Big(\sum_i(n_i-1)D_i\Big).
\end{align*}

Now by definition,
\[S_{(n_i-1)\Del_i}=\{\chi\in G^* : 1\le\tfrac{n_i}{d}\chi(\Fie(\Del_i))\le n_i-1\}=\{\chi\in G^*:\chi(\Fie(\Del_i))\ne0\}.\]
Thus in the decomposition of $\fie_*\Oh_X\big(\sum_i(n_i-1)D_i\big)$ given by Theorem \ref{thm!pushforward}, the summand $\eL_{\chi}^{-1}$ is twisted by $\sum_{j\in J}\Del_j$, where $J$ is the set of indices $j$ with $\chi(\Fie(\Del_j))\ne 0$. Then by \eqref{eqn!pardini},
\[\eL_{\chi}^{-1}\Big(\sum_{i\in J}\Del_i\Big)=\sum_i(1-\tfrac1d)\chi(\Fie(\Del_i))\Del_i=\eL_{\chi^{-1}},\]
where the last equality is because $\chi^{-1}(g)=-\chi(g)=d-\chi(g)$ for any $g$ in $G$. Thus we obtain
\[\fie_*\Big(\Oh_X\Big(\sum_i(n_i-1)D_i\Big)\Big)=\bigoplus_{\chi\in G^*}\eL_{\chi^{-1}},\]
and the Corollary follows. $\qed$
\end{pf}
\subsubsection{Example \ref{exa!Campedelli} continued}\label{exa!Campedelli-cont}
We resume our discussion of the Campedelli surface. The fundamental group of $X$ is $(\ZZ/2)^3$, and so the maximal abelian cover $\pi\colon A\to X$ is a $(\ZZ/2)^6$-cover $\psi\colon A\to\PP^2$ branched over $\Del$. Choose generators $g_1,\dots,g_6$ of $(\ZZ/2)^6$. As promised in Example \ref{exa!Campedelli}, we now fix $\Fie$ and $\Psi$:
\[\renewcommand{\arraystretch}{1.2}\begin{array}{c||c|c|c|c|c|c|c}
\Del_i&\Del_1&\Del_2&\Del_3&\Del_4&\Del_5&\Del_6&\Del_7\\
\hline\hline
\Phi(\Del_i)&g_1&g_2&g_3&g_1+g_2&g_1+g_3&g_2+g_3&g_1+g_2+g_3\\
\Psi(\Del_i)-\Phi(\Del_i)&0&0&0&g_4&g_5&g_6&g_4+g_5+g_6
\end{array}\]
For clarity, the table displays the difference between $\Psi(\Del_i)$ and $\Fie(\Del_i)$. In order that $A$ be the maximal abelian cover, $\Psi$ is defined so that each $\Psi(\Del_i)$ generates a distinct summand of $(\ZZ/2)^6$, excepting $\Psi(\Del_7)$, which is chosen so that $\sum_i\Psi(\Del_i)=0$. This last equality is induced by the relation $\sum_i\Del_i=0$ in $H_1(\PP^2-\Del,\ZZ)$.

The torsion group $\Tors X$ is generated by $g_4^*$, $g_5^*$, $g_6^*$. As an illustration of Theorem \ref{thm!pushforward}, we calculate $\fie_*\Oh_X(D_1)\otimes\tau$, where $\tau$ is the torsion line bundle on $X$ associated to $g_4^*$. Suppose $\fie_*\Oh_X(D_1)\otimes\tau=\bigoplus_{\chi\in G^*}\eM_{\chi}$, where $\eM_{\chi}$ are the line bundles to be calculated. In the table below, we collect the data relevant to Theorem \ref{thm!pushforward}.
\[\renewcommand{\arraystretch}{1.1}
\begin{array}{c|c|c|c|c}
\chi&\eL_{\chi+\tau}^{-1}&(\chi+\tau)\circ\Psi(D_1)&\text{Twist by }\Del_1?&\eM_{\chi}\\
\hline\hline
(0,0,0)&\Oh_{\PP^2}(-1)&0&\text{No}&\Oh_{\PP^2}(-1)\\
(1,0,0)&\Oh_{\PP^2}(-1)&1&\text{Yes}&\Oh_{\PP^2}\\
(0,1,0)&\Oh_{\PP^2}(-1)&0&\text{No}&\Oh_{\PP^2}(-1)\\
(0,0,1)&\Oh_{\PP^2}(-2)&0&\text{No}&\Oh_{\PP^2}(-2)\\
(1,1,0)&\Oh_{\PP^2}(-3)&1&\text{Yes}&\Oh_{\PP^2}(-2)\\
(1,0,1)&\Oh_{\PP^2}(-2)&1&\text{Yes}&\Oh_{\PP^2}(-1)\\
(0,1,1)&\Oh_{\PP^2}(-2)&0&\text{No}&\Oh_{\PP^2}(-2)\\
(1,1,1)&\Oh_{\PP^2}(-2)&1&\text{Yes}&\Oh_{\PP^2}(-1)
\end{array}\]

Summing the last column of the table, we get 
\[\fie_*\Oh_X(D_1)\otimes\tau=\Oh_{\PP^2}\oplus4\Oh_{\PP^2}(-1)\oplus3\Oh_{\PP^2}(-2).\]
In particular, we see that the linear system on $X$ associated to the line bundle $\Oh_X(D_1)\otimes\tau$ contains a single effective divisor.

\section{Exceptional collections of line bundles on surfaces}\label{sec!surfaces}
\subsection{Overview and definitions}\label{sec!strategy}
We outline our method for producing exceptional collections, starting with some definitions and fundamental observations. A good reference for semi-orthogonal decompositions is \cite{KuzICM}, and Proposition \ref{prop!K-theory} is proved in \cite{GO}.
\begin{dfn}\label{def!exceptional} An object $E$ in $D^b(X)$ is called \emph{exceptional} if
\[{\Ext}^k(E,E)=\left\{\begin{array}{ll}\CC&\text{ if }k=0,\\0&\text{ otherwise.}\end{array}\right.\]
An \emph{exceptional collection} $\EE\subset D^b(X)$ is a sequence of exceptional objects $\EE=(E_0,\dots,E_n)$ such that if $0\le i< j\le n$ then
\[\Ext^k(E_j,E_i)=0\text{ for all }k.\]
\end{dfn}
\begin{rmk} Some authors prefer the term exceptional sequence rather than exceptional collection.
\end{rmk}
It follows from Definition \ref{def!exceptional} that a line bundle on a surface is exceptional if and only if $p_g=q=0$. Moreover, if $\EE$ is an exceptional collection of line bundles, and $L$ is any line bundle, then $\EE\otimes L=(E_0\otimes L,\dots,E_n\otimes L)$ is again an exceptional collection, so we always normalise $\EE$ so that $E_0=\Oh_X$.

Let $\eE=\Span{\EE}$ denote the smallest full triangulated subcategory of $D^b(X)$ containing all objects in $\EE$. Then $\eE$ is an admissible subcategory of $D^b(X)$, and so we have a \emph{semiorthogonal decomposition}
\[D^b(X)={\Span{\eE,\Ay}},\]
where $\Ay$ is the left orthogonal to $\eE$. That is, $\Ay$ consists of all objects $F$ in $D^b(X)$ such that $\Ext^k(F,E)=0$ for all $k$ and for all $E$ in $\eE$. We say that the exceptional collection $\EE$ is \emph{full} if $D^b(X)=\eE$. The $K$-theory is additive across semiorthogonal decompositions:
\begin{prop}\label{prop!K-theory} If $D^b(X)=\Span{\Ay,\mathcal{B}}$ is a semiorthogonal decomposition, then
\[K_0(X)=K_0(\Ay)\oplus K_0(\mathcal{B}).\]
\end{prop}
% see Gorchinsky--Orlov or Alexeev--Orlov stated as a remark

Moreover, if $\EE$ is an exceptional collection of length $n$, then $K_0(\eE)=\ZZ^n$. Thus if $K_0(X)$ is not free, then $X$ can never have a full exceptional collection. The maximal length of an exceptional collection on $X$ is less than or equal to the rank of $K(X)$.
\subsubsection{Exceptional collections on del Pezzo surfaces}
Let $Y$ be the blow up of $\PP^2$ in $n$ points, and write $H$ for the pullback of the hyperplane section, $\Ebar_i$ for the $i$th exceptional curve. Then by work of Kuleshov and Orlov \cite{O}, \cite{KO} there is an exceptional collection of sheaves on $Y$
\[\Oh_{\Ebar_1}(-1),\dots,\Oh_{\Ebar_n}(-1),\Oh_Y,\Oh_Y(H),\Oh_Y(2H).\]
Note that the blown up points do not need to be in general position, and can even be infinitely near. 
We prefer an exceptional collection of line bundles on $Y$, so we mutate past $\Oh_Y$ to get
\begin{equation}\label{eqn!excdP}
\Oh_Y,\ \Oh_Y(\Ebar_1),\dots,\Oh_Y(\Ebar_n),\ \Oh_Y(H),\ \Oh_Y(2H).
\end{equation}
In fact, we only use the numerical properties of a given exceptional collection of line bundles on $Y$. Choose a basis $e_0,\dots,e_n$ for the lattice $\Pic Y\isom\ZZ^{1,n}$ with intersection form $\diag(1,-1^n)$. Then we write equation \eqref{eqn!excdP} numerically as
\[0,\ e_1,\dots,e_n,\ e_0,\ 2e_0.\]

\subsubsection{From del Pezzo to general type}
Let $X$ be a surface of general type with $p_g=0$ which admits an abelian cover $\fie\colon X\to Y$ of a del Pezzo surface $Y$ with $K_Y^2=K_X^2$. In addition, we suppose that the maximal abelian cover $A\to X\to Y$ is also Galois. Otherwise choose a maximal subgroup $T\subset\Tors X$ for which the  associated cover is Galois, and replace $A$, as in Section \ref{sec!prelim}. The branch divisor is $\Del=\sum_i\Del_i$ and we assume that $\Del$ is sufficiently reducible so that
\begin{enumerate}
\item[(A1)] $\Pic Y$ is generated by integral linear combinations of $\Del_i$.
\end{enumerate}
Now the Picard lattices of $X$ and $Y$ are isomorphic. Thus if $G$ is not too complicated, e.g.~of the form $\ZZ/p\times\ZZ/q$, we might hope to have:
\begin{enumerate}
\item[(A2)] The reduced pullbacks $D_i$ of $\Del_i$ (see Definition \ref{def!redpull}) generate $\Pic X/{\Tors X}$.
\end{enumerate}
In very good cases, reduced pullback actually induces an isometry of lattices
\begin{enumerate}
\item[(A3)] $f\colon\Pic Y\to\Pic X/{\Tors X}$, such that $f(K_Y)=-K_X$ modulo $\Tors X$.
\end{enumerate}
We say that a surface satisfies assumption (A) if (A1), (A2) and (A3) hold. These conditions are quite strong, and are not strictly necessary for our methods. For example, we could replace (A3) with an isometry of lattices from the abstract lattice $\ZZ^{1,n}$ to $\Pic X/{\Tors X}$.
%See Section \ref{sec!Keum-Naie} for such an example.

\begin{dfn} A sequence $\EE=(E_0,\dots,E_n)$ of line bundles on $X$ is called \emph{numerically exceptional} if $\chi(E_j,E_i)=0$ whenever $0\le i< j\le n$.
\end{dfn}

Assume $X$ satisfies (A), and let $(\Lambda_i)=(\Lambda_0,\dots,\Lambda_n)$ be an exceptional collection on $Y$. Now define $(L_i)=(L_0,\dots,L_n)$ by $L_i=f(\Lambda_i)^{-1}$. A calculation with the Riemann--Roch formula shows that $(L_i)$ is a numerically exceptional collection on $X$. This is explained in \cite{AO}.

Given a numerically exceptional collection $(L_i)$ of line bundles on $X$, the remaining obstacle is to determine whether $(L_i)$ is genuinely exceptional rather than just numerically so. Indeed, most numerically exceptional collections on $X$ are not exceptional. The standard trick (see \cite{BBS}) is to choose torsion line bundles $\tau_i$ in such a way that the twisted sequence $(L_i\otimes\tau_i)$ is an exceptional collection. We examine these choices of $\tau_i$ more carefully in what follows.

\subsubsection{Acyclic line bundles}
We discuss acyclic line bundles following \cite{GS}.
\begin{dfn} Let $L$ be a line bundle on $X$. If $H^i(X,L)=0$ for all $i$, then we call $L$ an \emph{acyclic} line bundle. We define the acyclic set associated to $L$ to be
\[\Ay(L)=\left\{\tau\in\Tors X : L\otimes\tau\text{ is acyclic}\right\}.\]
We call $L$ \emph{numerically acyclic} if $\chi(X,L)=0$. Clearly, an acyclic line bundle must be numerically acyclic.
\end{dfn}
\begin{rmk}\rm In the notation of \cite{GS}, $\tau=-\chi$.
\end{rmk}
\begin{lemma}[\cite{GS}, Lemma 3.4]\label{lem!acyclic} A numerically exceptional collection $L_0=\Oh_X,L_1\otimes\tau_1,\dots,$ $L_n\otimes\tau_n$ on $X$ is exceptional if and only if
\begin{equation}
\begin{split}\label{eqn!torsion-conditions}
-\tau_i&\in\Ay(L_i^{-1})\text{ for all }i, \text{ and}\\
\tau_i-\tau_j&\in\Ay(L_j^{-1}\otimes L_i)\text{ for all }j>i.
\end{split}
\end{equation}
\end{lemma}

Thus to enumerate all exceptional collections on $X$ of a particular numerical type, it suffices to calculate the relevant acyclic sets, and systematically test the above conditions \eqref{eqn!torsion-conditions} on all possible combinations of $\tau_i$.

\subsubsection{Calculating cohomology of line bundles}\label{sec!cohom}
Given a torsion twist $L\otimes\tau$, Theorem \ref{thm!pushforward} gives a decomposition
\[\fie_*(L\otimes\tau)=\bigoplus_{\chi\in G^*}\eM_\chi,\]
for some line bundles $\eM_\chi$ on $Y$, which may be computed explicitly. Since $\fie$ is finite, we have
\[h^p(L\otimes\tau)=\sum_{\chi\in G^*}h^p(\eM_\chi)\]
for all $p$.

Thus $L\otimes\tau$ is acyclic if and only if each summand $\eM_\chi$ is acyclic on $Y$. Now if $\chi(Y,\eM_\chi)=0$ and $h^0(\eM_\chi)=h^2(\eM_\chi)=0$, we see that $h^1(\eM_\chi)=0$. Thus by Serre duality and the Riemann--Roch theorem, we are reduced to calculating Euler characteristics and determining effectivity for (lots of) divisor classes on the del Pezzo surface $Y$.

\subsubsection{Coordinates on $\Pic X/{\Tors X}$}
Under assumption (A), we make the following definition.
\begin{dfn}\label{def!coordinates} Choose a basis $B_1,\dots,B_n$ for $\Pic X/{\Tors X}$ consisting of linear combinations of reduced pullbacks. Then any line bundle $L$ on $X$ may be written uniquely as
\[L=\Oh_X(d_1,\dots,d_n)\otimes\tau\]
so that $L=\Oh_X\big(\sum_{i=1}^nd_iB_i\big)\otimes\tau$. We call $d$ (respectively $\tau$) the multidegree (resp.~torsion twist) of $L$ with respect to the chosen basis.
\end{dfn}
The torsion twist associated to any line bundle on $X$ may be calculated using Theorem \ref{thm!pushforward} and the following immediate lemma. See Lemma \ref{lem!Kulikov-torsion} for an example.
\begin{lemma}\label{lem!cohomology-torsion} If $\tau$ is a torsion line bundle, then $h^0(\tau)\ne 0$ implies $\tau=0$.
\end{lemma}
\begin{rmk}\rm Definition \ref{def!coordinates} fixes a basis for $\Pic Y=\ZZ^{1,9-K^2}$ via the isometry with $\Pic X/{\Tors X}$. This basis corresponds to a geometric marking on the del Pezzo surface $Y$, and the multidegree $d$ of $L$ is just the image of $L$ in $\Pic Y$ under the isometry. In fixing our basis, we break some of the symmetry of the coordinates. This is necessary in order to use the computer to search for exceptional collections. We can recover the symmetry later using the Weyl group action (see Section \ref{sec!group-actions}).
\end{rmk}

\subsubsection{Determining effectivity of divisor classes}\label{sec!effective}

For each fake del Pezzo surface, we have the following theorem.

\begin{thm} Suppose $X$ is a fake del Pezzo surface satisfying assumption $(A)$ and with $T=\Tors X$. Let $\Eff$ denote the semigroup generated by the reduced pullbacks $D_i$ of the irreducible branch components $\Del_i$, and pullbacks of the other $(-1)$- and $(-2)$-curves on $Y$. Then $\Eff$ is the semigroup of all effective divisors on $X$.
\end{thm}
We prove this theorem for the secondary nodal Burniat surface with $K^2=4$ in Section \ref{sec!burniat-secondary} (cf.~\cite{Alexeev} for the Burniat surface with $K^2=6$). The other fake del Pezzo surfaces work in the same way, see \cite{computer-code}.

Moreover, $\Eff$ is graded by multidegree, and we define a homomorphism 
\[t\colon\Eff\to\Tors X\]
sending $D_i$ to its torsion twist under Definition \ref{def!coordinates}. The image under $t$ of the graded summand $\Eff_d$ of multidegree $d$ is the set of torsion twists $\tau$ for which $\Oh_X(\sum d_iB_i)\otimes\tau$ is effective.

\subsubsection{Group actions on the set of exceptional collections}\label{sec!group-actions}

We consider a dihedral group action and the Weyl group action on the set of exceptional collections on $X$. Mutations are not considered systematically in this article, since a mutation of a line bundle need not be a line bundle.

Let $\EE=(E_1,\dots,E_n)$ be an exceptional collection of line bundles on $X$. If we normalise the first line bundle of any exceptional collection to be $\Oh_X$, then there is an obvious dihedral group action on the set of exceptional collections of length $n$ on $X$, generated by $\EE\mapsto(E_2,\dots,E_n,E_1(-K_X))$ and $\EE\mapsto\EE^{-1}=(E_n^{-1},\dots,E_1^{-1})$.

The Weyl group of $\Pic Y$ is generated by reflections in $(-2)$-classes. That is, suppose $\al$ is a class in $\Pic Y$ with $K_Y\cdot\al=0$ and $\al^2=-2$. Then
\[r_{\al}\colon L\mapsto L+(L\cdot\al)\al\]
is a reflection on $\Pic Y$ which fixes $K_Y$. Any reflection sends an exceptional collection on $Y$ to another exceptional collection. Thus the Weyl group action on numerical exceptional collections on $Y$ induces an action on numerical exceptional collections on $X$ under assumption (A). This action accounts for the choices made in giving $Y$ a geometric marking (see Definition \ref{def!coordinates}).

\subsection{The Kulikov surface with $K^2=6$}\label{sec!Kulikov}
\begin{figure}[ht]
\begin{picture}(100,70)(-120,5)
\put(10,10){\line(1,0){80}}
\put(90,10){\line(-2,3){40}}
\put(10,10){\line(2,3){40}}
\put(10,10){\line(5,3){65}}
\put(90,10){\line(-5,2){70}}
\put(61,4){\line(-1,6){11}}
\put(10,10){\circle*{5}}
\put(90,10){\circle*{5}}
\put(50,70){\circle*{5}}
\put(25,1){$\Del_1$}
\put(80,27){$\Del_2$}
\put(23,55){$\Del_3$}
\put(60,0){$\Del_4$}
\put(76,47){$\Del_5$}
\put(5,35){$\Del_6$}
\put(47,74){$P_1$}
\put(-5,5){$P_2$}
\put(92,5){$P_3$}
\end{picture}
\caption{The Kulikov configuration}\label{fig!kulikov}
\end{figure}
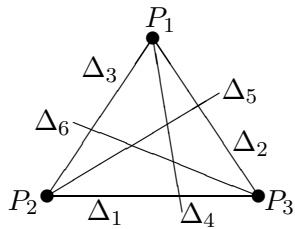

For details on the Kulikov surface (first described in \cite{Ku}), its torsion group and moduli space, see \cite{CC}. The Kulikov surface $X$ is a $(\ZZ/3)^2$-cover of the del Pezzo surface $Y$ of degree $6$. Figure \ref{fig!kulikov} shows the associated cover of $\PP^2$ branched over six lines in special position. The configuration has just one free parameter, and in fact, the Kulikov surfaces form a $1$-dimensional, irreducible, connected component of the moduli space of surfaces of general type with $p_g=0$ and $K^2=6$.

To obtain a nonsingular cover, we blow up the plane at three points $P_1,P_2,P_3$, giving a $(\ZZ/3)^2$-cover of a del Pezzo surface of degree $6$. The exceptional curves are denoted $\Ebar_i$. By results of \cite{CC}, the torsion group $\Tors X$ is isomorphic to $(\ZZ/3)^3$, so the maximal abelian cover $\psi\colon A\to Y$ has group $\Gtilde\isom(\ZZ/3)^5$. Let $g_i$ generate $\Gtilde$, and write $g_i^*$ for the dual generators of $\Gtilde^*$. As explained in Section \ref{sec!prelim}, the covers are determined by $\Fie\colon H_1(\PP^2-\Del,\ZZ)\to G$ and $\Psi\colon H_1(\PP^2-\Del,\ZZ)\to \Gtilde$, which are defined in the table below.
\[\renewcommand{\arraystretch}{1.3}\begin{array}{c||c|c|c|c|c|c}
D&\Del_1&\Del_2&\Del_3&\Del_4&\Del_5&\Del_6\\
\hline\hline
\Fie(D)&g_1&g_1&g_1&g_2&g_1+g_2&2g_1+g_2\\
\Psi(D)-\Fie(D)&0&g_3&2g_3+g_4&2g_4&g_5&2g_5
\end{array}\]
 The images of the exceptional curves $\Ebar_i$ under $\Fie$ and $\Psi$ are computed using formula \eqref{eqn!exceptional-curve}:
\[\Fie(\Ebar_1)=2g_1+g_2,\ \Fie(\Ebar_2)=g_2,\ \Fie(\Ebar_3)=g_1+g_2, \text{ etc.}\]
\begin{lemma}\label{lem!Kulikov-picard}
The Kulikov surface satisfies assumptions (A1) and (A2). That is, the free part of $\Pic X$ is generated by the reduced pullbacks of $\Del_1+\Ebar_2+\Ebar_3$, $\Ebar_1$, $\Ebar_2$, $\Ebar_3$, and the intersection pairing $\diag(1,-1,-1,-1)$ is inherited from $Y$.
\end{lemma}
\begin{pf} Define $e_0=D_1+E_2+E_3$, $e_1=E_1$, $e_2=E_2$, $e_3=E_3$ in $\Pic X$. These are integral divisors, since they are reduced pullbacks, and the intersection pairing is $\diag(1,-1,-1,-1)$, which is unimodular. For example, by definition of reduced pullback, $3e_0=\fie^*(\Del_1+\Ebar_2+\Ebar_3)$, and so
\[(3e_0)^2=\fie^*(\Del_1+\Ebar_2+\Ebar_3)^2=9\cdot1,\]
or $e_0^2=1$. Hence we have an isomorphism of lattices.$\qed$
\end{pf}
Using the basis chosen in this lemma, we compute the coordinates (Definition \ref{def!coordinates}) of the reduced pullback $D_i$ of each irreducible branch component $\Del_i$.
\begin{lemma}\label{lem!Kulikov-torsion} We have
\begin{align*}
\Oh_X(D_1)&=\Oh_X(1,0,-1,-1),& \Oh_X(D_4)&=\Oh_X(1,-1,0,0)[2,1,2],\\
\Oh_X(D_2)&=\Oh_X(1,-1,0,-1)[1,0,2],&\Oh_X(D_5)&=\Oh_X(1,0,-1,0)[2,1,0],\\
\Oh_X(D_3)&=\Oh_X(1,-1,-1,0)[2,0,2],&\Oh_X(D_6)&=\Oh_X(1,0,0,-1)[2,1,1],
\end{align*}
where $[a,b,c]$ in $(\ZZ/3)^3$ denotes a torsion line bundle on $X$.
\end{lemma}
\begin{pf}
We prove that $\Oh_X(D_2)=\Oh_X(1,-1,0,-1)[1,0,2]$. The other cases are similar. It is clear that $\Del_2\lin\Del_1-\Ebar_1+\Ebar_2$ on $Y$, so the multidegree is correct. It remains to check the torsion twist, by showing that $\eF=\Oh_X(D_2-D_1+E_1-E_2-\tau)$ has a global section when $\tau=[1,0,2]$. Then by Lemma \ref{lem!cohomology-torsion}, we have the desired equality.

The pushforward $\fie_*\eF$ splits into a direct sum of line bundles $\bigoplus\eM_{\chi}$, one for each character $\chi=(a,b)$ in $G^*$. The following table collects the data required to calculate each $\eM_{\chi}$ via Theorem \ref{thm!pushforward}. The second column is calculated using equation \eqref{eqn!pardini}, and the next four columns evaluate $\chi+\tau$ on each $\Psi(\Gam)$, where $\Gam$ is any one of $\Del_1$, $\Del_2$, $\Ebar_1$ and $\Ebar_2$. The final column is explained below.
\[\renewcommand{\arraystretch}{1.3}
\begin{array}{c|c|c|c|c|c|c}
&&\multicolumn{4}{c|}{(\chi+\tau)\circ\Psi(\Gam)}& \\
\cline{3-6}
\chi&\eL_{\chi+\tau}^{-1}&\Del_1&\Del_2&\Ebar_1&\Ebar_2&\eM_\chi\\
\hline\hline
(0,0)&\Oh_Y(-2,1,1,0)&0&1&0&1&\Oh_Y(-3,1,2,1)\\
(1,0)&\Oh_Y(-1,0,0,1)&1&2&2&1&\Oh_Y\\
(0,1)&\Oh_Y(-2,1,0,1)&0&1&1&2&\Oh_Y(-3,1,1,2)\\
(2,0)&\Oh_Y(-2,0,1,1)&2&0&1&1&\Oh_Y(-2,0,1,1)\\
(1,1)&\Oh_Y(-2,1,0,1)&1&2&0&2&\Oh_Y(-1,0,0,0)\\
(0,2)&\Oh_Y(-2,1,1,0)&0&1&2&0&\Oh_Y(-3,2,1,1)\\
(2,1)&\Oh_Y(-2,0,1,0)&2&0&2&2&\Oh_Y(-2,1,1,0)\\
(1,2)& \Oh_Y(-3,1,1,1)&1&2&1&0&\Oh_Y(-2,0,0,0)\\
(2,2)& \Oh_Y(-2,1,1,1)&2&0&0&0&\Oh_Y(-2,1,0,1)
\end{array}\]

Now by the projection formula (cf.~Remark \ref{rmk!projformula}),
\[\fie_*\eF=\fie_*\Oh_X(2D_1+D_2+E_1+2E_2-\tau)\otimes\Oh_Y(-\Del_1-\Ebar_2).\]
So according to Theorem \ref{thm!pushforward} and the remark following it, each $\eM_\chi$ is a twist of $\eL_{\chi+\tau}^{-1}(-\Del_1-\Ebar_2)$ by a certain combination of $\Del_1$, $\Del_2$, $\Ebar_1$ and $\Ebar_2$. By Definition \ref{def!S_i}, the rules governing the twists are:
\begin{align*}
\text{twist by }\Del_1&\iff(\chi+\tau)\circ\Psi(\Del_1)=1\text{ or }2\\
\text{twist by }\Del_2&\iff(\chi+\tau)\circ\Psi(\Del_2)=2\\
\text{twist by }\Ebar_1&\iff(\chi+\tau)\circ\Psi(\Ebar_1)=2\\
\text{twist by }\Ebar_2&\iff(\chi+\tau)\circ\Psi(\Ebar_2)=1\text{ or }2.
\end{align*}

Thus $\fie_*\eF$ is given by the direct sum of the line bundles $\eM_{\chi}$ listed in the final column. Note that $\eM_{(1,0)}=\Oh_Y$, so $h^0(\fie_*\eF)=1$. Hence $D_2-D_1+E_1-E_2-\tau\lin 0$.
$\qed$
\end{pf}
\begin{cor} By formula \eqref{eqn!canonical}, we have
\[\Oh_X(K_X)=\Oh_X(3,-1,-1,-1)[0,0,2].\]
Thus the Kulikov surface satisfies (A3).
\end{cor}
\begin{pf} The multidegree is clear by \eqref{eqn!canonical}, but the torsion twist requires some care. Since $K_X$ is the pullback of an integral divisor on $Y$, it should be torsion-neutral with respect to our coordinate system on $\Pic X$. Thus by Lemma \ref{lem!Kulikov-torsion}, we see that the required twist is $[0,0,2]$.$\qed$
\end{pf}
\begin{thm}\label{thm!Kulikov-effective} The semigroup $\Eff$ of effective divisors on the Kulikov surface is generated by the nine reduced pullbacks of components of the branch divisor $D_1,\dots,D_6,E_1,E_2,E_3.$ $\qed$
\end{thm}
This Theorem is proved using an easier variant of the proof of Theorem \ref{thm!Burniat-effective}. The situation here is easier, because all of the $(-1)$-curves on $Y$ are branch divisors, and there are no $(-2)$-curves.

Thus we have a homomorphism of semigroups $t\colon\Eff\to\Tors X$, which sends an effective divisor to its associated torsion twist (see Lemma \ref{lem!Kulikov-picard}), under the choice of basis (from Lemma \ref{lem!Kulikov-torsion}).

\subsubsection{Acyclic line bundles on the Kulikov surface}

Let us start with the following numerical exceptional collection on $Y$:
\[\Lambda\colon 0,\ e_0-e_1,\ e_0-e_2,\ e_0-e_3,\ 2e_0-\textstyle{\sum}_{i=1}^3 e_i,\ e_0.\]
Given assumptions (A), we see that $\La$ corresponds to the following numerically exceptional sequence of line bundles on $X$:
\begin{equation}\label{eqn!Kulikov-num-exc}
\begin{split}
&L_0=\Oh_X,\ L_1=\Oh_X(-1,1,0,0),L_2=\Oh_X(-1,0,1,0),\\
L_3=&\Oh_X(-1,0,0,1),\ L_4=\Oh_X(-2,1,1,1),\ L_5=\Oh_X(-1,0,0,0).
\end{split}
\end{equation}

We find all collections of torsion twists $L_i\otimes\tau_i$ which are exceptional collections on $X$. The first step is to find the acyclic sets associated to the various $L_j^{-1}\otimes L_i$.
\begin{prop}\label{prop!Kulikov-acyclic} The acyclic sets $\Ay(L_j^{-1}\otimes L_i)$ for $j>i\ge0$ are listed in Appendix \ref{app!Kulikov}.
\end{prop}
\begin{pf1}
By Theorem \ref{thm!Kulikov-effective}, it is an easy exercise to check each entry in the table. As an illustration, we calculate $\Ay(L_1^{-1})$. The effective divisors on $X$ of multidegree $(1,-1,0,0)$ are $D_2+E_3$, $D_3+E_2$, $D_4$.
Thus applying the homomorphism $t$ to each of these effective divisors, we see that $[1,0,2]$, $[2,0,2]$, $[2,1,2]$ do not appear in $\Ay(L_1^{-1})$.
Next we consider degree two cohomology via Serre duality. The effective divisors of multidegree $(2,0,-1,-1)$ are
\begin{gather*}
2D_1+E_2+E_3,D_1+D_2+E_1+E_3,D_1+D_3+E_1+E_2,\\
D_2+D_3+2E_1,D_1+D_4+E_1,D_1+D_5+E_2,D_1+D_6+E_3,\\
D_2+D_5+E_1,D_3+D_6+E_1,D_5+D_6.
\end{gather*}
Again, applying $t$ we find that $[0,0,2]$, $[2,0,0]$, $[1,0,0]$, $[0,0,1]$, $[1,2,0]$, $[1,2,2]$, $[1,2,1]$, $[0,2,0]$, $[2,2,2]$, $[2,1,1]$ can not appear in $\Ay(L_1^{-1})$. The acyclic set is made up of those elements of $\Tors X$ which do not appear in either of the two lists above.$\qed$
\end{pf1}
\begin{pf2}
As a sanity check, an alternative proof is to use Theorem \ref{thm!pushforward} repeatedly, to calculate the cohomology of all possible torsion twists of $L_1$.
$\qed$
\end{pf2}
Both methods are implemented in our computer script \cite{computer-code}.

\subsubsection{Exceptional collections on the Kulikov surface}
We now find all exceptional collections on $X$ which are numerically of the form \eqref{eqn!Kulikov-num-exc}. Lemma \ref{lem!acyclic} reduces us to a simple search, which can be done systematically \cite{computer-code}.
\begin{thm}\label{thm!Kulikov-exceptional} The surface $X$ has nine exceptional collections $L_0=\Oh_X$, $L_1\otimes\tau_1,\dots,L_5\otimes\tau_5$ which are numerically of the form \eqref{eqn!Kulikov-num-exc}. They are given in Table \ref{tab!Kulikov-exceptional} below. Each row lists the required torsion twists $\tau_i$ for $i=1,\dots,5$ as elements of $(\ZZ/3)^3$.
\end{thm}
\begin{table}[h]
\[\renewcommand{\arraystretch}{1.3}
\begin{array}{c|ccccc}
&\tau_1&\tau_2&\tau_3&\tau_4&\tau_5\\
\hline\hline
1&[ 0, 0, 0 ]& [ 0, 2, 2 ]& [ 2, 2, 1 ]& [ 2, 2, 1 ]& [ 0, 0, 1 ]\\
2&[ 2, 2, 0 ]& [ 2, 1, 2 ]& [ 0, 0, 1 ]& [ 1, 1, 1 ]& [ 2, 2, 1 ]\\
3&[ 2, 2, 1 ]& [ 2, 1, 2 ]& [ 0, 0, 1 ]& [ 1, 1, 1 ]& [ 2, 0, 2 ]\\
4&[ 2, 2, 0 ]& [ 2, 0, 1 ]& [ 0, 2, 0 ]& [ 2, 2, 1 ]& [ 2, 1, 2 ]\\
5&[ 1, 1, 0 ]& [ 1, 0, 2 ]& [ 2, 2, 0 ]& [ 1, 1, 1 ]& [ 2, 2, 1 ]\\
6&[ 1, 1, 0 ]& [ 1, 0, 2 ]& [ 0, 0, 1 ]& [ 1, 1, 1 ]& [ 2, 2, 1 ]\\
7&[ 1, 1, 0 ]& [ 1, 0, 2 ]& [ 2, 2, 1 ]& [ 1, 1, 1 ]& [ 0, 0, 1 ]\\
8&[ 2, 0, 2 ]& [ 2, 2, 0 ]& [ 0, 1, 2 ]& [ 1, 1, 1 ]& [ 2, 2, 1 ]\\
9&[ 2, 0, 2 ]& [ 2, 2, 1 ]& [ 0, 1, 2 ]& [ 1, 1, 1 ]& [ 1, 0, 2 ]
\end{array}
\]
\caption{Exceptional collections on the Kulikov surface}\label{tab!Kulikov-exceptional}
\end{table}
\begin{rmk}\rm
\begin{enumerate}
\item The precise number of exceptional collections is not important. Rather, the fact that we have definitively enumerated all exceptional collections of numerical type $\La$, means that we can sift through the list to find one with the most desirable properties.
\item Let $\La'$ be any translation of $\La$ under the Weyl group action of $A_1\times A_2$ on $\Pic Y$. Then $\La'$ is another numerical exceptional collection on $X$ (see Section \ref{sec!group-actions}), so we may enumerate exceptional collections on $X$ of numerical type $\La'$. For the Kulikov surface, each element of the orbit corresponds to either 9, 14, 18 or 24 exceptional collections on $X$. Thus, the Weyl group action does not ``lift'' to $X$ in a way which is compatible with the covering $X\to Y$. On occasion, this incompatibility is used to our advantage (see \cite{computer-code}). We return to these exceptional collections in Section \ref{sec!heights}.
\end{enumerate}
\end{rmk}

\section{Heights of exceptional collections}\label{sec!heights}
Let $X$ be a surface of general type with $p_g=q=0$, $\Tors X\ne 0$ with an exceptional collection of line bundles $\EE=(E_0,\dots,E_{n-1})$. Write $\eE$ for the smallest full triangulated subcategory of $D^b(X)$ containing $\EE$. In this section we calculate some invariants of $\EE$. The invariants we consider are essentially determined by the derived category, but we must enhance the derived category in order to make computations. For completeness, we discuss some background first.

\subsection{Motivation from del Pezzo surfaces}\label{sec!dg-motivate}
Let $Y$ be a del Pezzo surface and let $\EE$ be a strong exceptional collection of line bundles on $Y$. Recall that $\EE$ is \emph{strong} if $\Ext^k(E_i,E_j)=0$ for all $i$, $j$ and for all $k>0$. We define the \emph{partial tilting bundle} of $\EE$ to be $T=\bigoplus_i E_i$. Then the derived endomorphism ring $\Ext^*(T,T)=\bigoplus_{i,j}\Hom(E_i,E_j)$ is an associative algebra, and we have an equivalence of categories $\eE\isom D^b(\text{mod-}\Ext^*(T,T))$ (see \cite{BK}).

From now on, we assume that $\EE$ is an exceptional collection on a fake del Pezzo surface $X$, so that we do not have the luxury of choosing a strong exceptional collection. Instead, we recover $\eE$ by studying the higher multiplications coming from the $A_\infty$-algebra structure on $\Ext^*(T,T)$.

\subsection{Digression on dg-categories}\label{sec!dgcat}
We sketch the construction of a \emph{differential graded} (or dg) enhancement $\De$ of $D^b(X)$. Objects in $\De$ are the same as those in $D^b(X)$, but morphisms $\Hom^\bullet_{\De}(F,G)$ form a chain complex, with differential $d$ of degree $+1$. Composition of maps $\Hom^\bullet_\De(F,G)\otimes\Hom^\bullet_\De(G,H)\to\Hom^\bullet_\De(F,H)$ is a morphism of complexes (the Leibniz rule), and for any object $F$ in $\De$, we require $d(\id_F)=0$. For a precise definition of $\Hom^\bullet_\De(F,G)$, one could use the \v{C}ech complex, and we refer to \cite{K} for details. The main point is that the cohomology of $\Hom^\bullet_{\De}(F,G)$ in degree $k$ is $\Ext^k_{D^b(X)}(F,G)$, so in particular, we have $H^0(\Hom^\bullet_\De(F,G))=\Hom_{D^b(X)}(F,G)$. 

\subsection{Hochschild homology}

We first compute some additive invariants, only making implicit use of the dg-structure. The Hochschild homology of $X$ is given by the Hochschild--Kostant--Rosenberg isomorphism
\[\HH_k(X)\isom\bigoplus_pH^{p+k}(X,\Omega^p_X),\]
so $\HH_0(X)=\CC^{12-K^2}$ and $\HH_k(X)=0$ in all other degrees. Moreover, Hochschild homology is additive over semiorthogonal decompositions.
\begin{thm}\cite{Kuz09} If $D^b(X)=\Span{\Ay,\mathcal{B}}$ is a semiorthogonal decomposition, then
\[\HH_k(X)=\HH_k(\Ay)\oplus\HH_k(\mathcal{B}).\]
\end{thm}
Assuming the Bloch conjecture on algebraic zero-cycles, we have
\[K_0(X)=\ZZ^{12-K^2}\oplus\Tors X,\]
and we note that $K$-theory is also additive over semiorthogonal decompositions (see Proposition \ref{prop!K-theory}).

Now for an exceptional collection of length $n$, $K_0(\eE)=\ZZ^n$ and 
\[\HH_k(\eE)=\left\{\begin{array}{ll}\CC^n&\text{ if }k=0\\0&\text{ otherwise.}\end{array}\right.\]
Thus the maximal length of $\EE$ is at most $12-K_X^2$, and such an exceptional sequence of maximal length effects a semiorthogonal decomposition $D^b(X)=\Span{\eE,\Ay}$ with nontrivial semiorthogonal complement $\Ay$. We say that $\Ay$ is a \emph{quasiphantom} category; by additivity, the Hochschild homology vanishes, but $K_0(\Ay)\supseteq\Tors X\ne0$, so $\Ay$ can not be trivial.

\subsection{Height}\label{sec!dgcat-height}
The Hochschild cohomology groups of $X$ may be computed via the other Hochschild--Kostant--Rosenberg isomorphism (cf.~\cite{Kuz09}):
\[\HH^k(X)=\bigoplus_{p+q=k}H^q(X,\La^pT_X).\]
Thus for a surface of general type with $p_g=0$, we have 
\begin{gather*}\HH^0(X)\isom H^0(\Oh_X)=\CC,\ \HH^1(X)=0,\ \HH^2(X)\isom H^1(T_X),\\
\HH^3(X)\isom H^2(T_X),\ \HH^4(X)\isom H^0(2K_X)=\CC^{1+K^2}.
\end{gather*}
Recall that the degree two (respectively three) Hochschild cohomology is the tangent space (resp.~obstruction space) to the formal deformations of a category \cite{KS}.

In principle, \cite{K} gives an algorithm for computing $\HH^*(\Ay)$ using a spectral sequence and the notion of height of an exceptional collection. Moreover, by \cite[Prop.~6.1]{K}, for an exceptional collection to be full, its height must vanish. Thus the height may be used to prove existence of phantom categories without reference to the $K$-theory. We outline the algorithm of \cite{K} below.

Given an exceptional collection $\EE$ on $X$, there is a long exact sequence (induced by a distinguished triangle)
\[\ldots\to\NHH^k(\EE,X)\to\HH^k(X)\to\HH^k(\Ay)\to\NHH^{k+1}(\EE,X)\to\ldots\]
where $\NHH(\EE,X)$ is the \emph{normal Hochschild cohomology} of the exceptional collection $\EE$. The normal Hochschild cohomology can be computed using a spectral sequence with first page
\begin{align*}
\textbf{E}^1_{-p,q}=\bigoplus_{\substack{0\le a_0<\dots< a_p\le n-1\\ k_0+\dots+k_p=q}}&\Ext^{k_0}(E_{a_0},E_{a_1})\otimes\cdots \\
&\dots\otimes\Ext^{k_{p-1}}(E_{a_{p-1}},E_{a_p})\otimes\Ext^{k_p}(E_{a_p},S^{-1}(E_{a_0})).
\end{align*}
The spectral sequence relies on the dg-structure on $\De$; the initial differentials $d'$ and $d''$ are induced by the diff\-erential on $\mathcal{D}$ and the composition map respectively, while the higher differentials are related to the $A_\infty$-algebra structure on $\Ext$-groups, (see Section \ref{sec!A-infinity}).

The existing examples of exceptional collections on surfaces of general type with $p_g=0$ suggest that $\NHH^k(\EE,X)$ vanishes for small $k$. Thus the \emph{height} $h(\EE)$ of an exceptional collection $\EE=(E_0,\dots,E_{n-1})$ is defined to be the smallest integer $m$ for which $\NHH^m(\EE,X)$ is nonzero. Alternatively, $m$ is the largest integer such that the canonical restriction morphism $\HH^k(X)\to\HH^k(\Ay)$ is an isomorphism for all $k\le m-2$ and injective for $k=m-1$.

\subsection{Pseudoheight}
The height may be rather difficult to compute in practice, requiring a careful analysis of the $\Ext$-groups of $\EE$ and the maps in the spectral sequence. The pseudoheight is easier to compute and sometimes gives a good lower bound for the height.
\begin{dfn}\label{def!pseudoheight} The \emph{pseudoheight} $ph(\EE)$ of an exceptional collection $\EE=(E_0,\dots,E_{n-1})$ is
\begin{align*}
ph(\EE)=\min_{0\le a_0<\dots<a_p\le n-1} \big(&e(E_{a_0},E_{a_1})+\cdots\\
&+e(E_{a_{p-1}},E_{a_p})+e(E_{a_p},E_{a_0}({-K_X}))-p+2\big),
\end{align*}
where $e(F,F')=\min\{i:\Ext^i(F,F')\ne0\}$.
\end{dfn}
The pseudoheight is just the total degree of the first nonzero term in the first page of the spectral sequence, where the shift by $2$ takes care of the Serre functor.

Consider the length $2n$ anticanonical extension of the sequence $\EE$ (see also Section \ref{sec!group-actions}):
\begin{equation}\label{eqn!antican-extend}E_0,\dots,E_{n-1},E_n=E_0(-K_X),\dots,E_{2n-1}=E_{n-1}(-K_X).
\end{equation}
If the $E_i$ are line bundles, then we have a numerical lower bound for the pseudoheight.
\begin{lemma}\cite[Lem.~4.10, Lem.~5.1]{K} If $K_X$ is ample and $E_i\cdot K_X\ge E_j\cdot K_X$ for all $i<j$ and for all $E_i$, $E_j$ in the anticanonically extended sequence \eqref{eqn!antican-extend}, then $ph(\EE)\ge3$.
\end{lemma}
The numerical conditions required by the Lemma are not particularly stringent. For example, all the exceptional collections we have exhibited on the Kulikov surface in Section \ref{sec!Kulikov} have pseudoheight at least $3$, even before we consider the $\Ext$-groups more carefully.
\begin{rmk}\rm If $L$ is a line bundle, then $\dim\Ext^k(L,L(-K_X))=h^{2-k}(2K_X)$ by Serre duality, which is the case $p=0$ in Definition \ref{def!pseudoheight}. Thus any exceptional collection of line bundles on a surface of general type with $p_g=0$ has pseudoheight at most $4$. Moreover, if $ph(\EE)=4$, then $h(\EE)=4$ by \cite{K}.
\end{rmk}
\subsection{The $A_\infty$-algebra of an exceptional collection}\label{sec!A-infinity}
Let $\EE=(E_0,\dots,E_{n-1})$ be an exceptional collection on $X$, and define $T=\oplus_{i=0}^{n-1}E_i$. Then $B=\Hom^\bullet_\De(T,T)$ is a differential graded algebra via the dg-structure on $\De$ (see Section \ref{sec!dgcat}). It can be difficult to compute the dg-algebra structure on $B$ directly, so we pass to the $A_\infty$-algebra $H^*B$.

We discuss $A_\infty$-algebras, referring to \cite{Keller} for details and further references. An \emph{$A_\infty$-algebra} is a graded vector space $A=\bigoplus_{p\in\ZZ}A^p$, together with graded multiplication maps $m_n\colon A^{\otimes n}\to A$ of degree $2-n$, for each $n\ge 1$. These multiplication maps satisfy an infinite sequence of relations, starting with
\begin{align*}
m_1m_1&=0,\\
m_1m_2&=m_2(m_1\otimes\id_A+\id_A\otimes m_1).
\end{align*}
These first two relations ensure that $m_1$ is a differential on $A$, satisfying the Leibniz rule with respect to $m_2$. The third relation is
\begin{align*}
m_2(&\id_A\otimes m_2-m_2\otimes\id_A)=\\
&m_1m_3+m_3(m_1\otimes\id_A\otimes\id_A+\id_A\otimes m_1\otimes\id_A+\id_A\otimes\id_A\otimes m_1),
\end{align*}
which shows that $m_2$ is not associative in general, but if $m_n=0$ for all $n\ge3$, then $A$ is an ordinary associative differential graded algebra.

In fact, by the above discussion, we can view $B$ as an $A_\infty$-algebra, with $m_1$ being the differential, $m_2$ the multiplication, and $m_n=0$ for $n\ge3$.
By a theorem of Kadeishvili (cf.~\cite{Keller}), the homology $H^*B=H^*(B,m_1)$ has a canonical $A_\infty$-algebra structure, for which $m_1=0$, $m_2$ is induced by the multiplication on $B$, and $H^*B$ and $B$ are quasi-isomorphic as $A_\infty$-algebras. This canonical $A_\infty$-structure is unique, and $H^*B$ is called a \emph{minimal model} for $B$. We say that $B$ is \emph{formal} if it has a minimal model $H^*B$ for which $m_n=0$ for all $n\ge3$, so that $H^*B$ is just an associative graded algebra.

The $A_\infty$-algebra of $\EE$ is
\[H^*B=\Ext^*(T,T)=\bigoplus_k\bigoplus_{0\le i,j\le n-1}\Ext^k(E_i,E_j),\]
and $m_2$ coincides with the Yoneda product on $\Ext$-groups. Clearly, if the exceptional collection $\EE$ consists of sheaves, then $H^*B$ has only three nontrivial graded summands, in degrees $0$, $1$ and $2$. Since $m_n$ has degree $2-n$, the summands of degree $0$ and $1$ are crucial in determining the $A_\infty$-algebra structure.

\subsubsection{Recovering $\eE$ from $H^*B$}
According to \cite{BK}, \cite{Keller94}, the subcategory $\eE$ of $\De$ generated by the exceptional collection $\EE$ is equivalent to the triangulated subcategory $\Perf(B)\subset D^b(\text{mod-}B)$ of perfect objects over the dg-algebra $B$. A perfect object is a differential graded $B$-module that is quasi-isomorphic to a bounded chain complex of projective and finitely generated modules. As mentioned above, it is preferable to consider the $A_\infty$-algebra $H^*B$ instead, noting that $\eE$ is in turn equivalent to the triangulated category of perfect $A_\infty$-modules over $H^*B$. If $B$ is formal, the equivalence reduces to $\eE\isom D^b(\text{mod-}H^*B)$, which should be compared with Section \ref{sec!dg-motivate}.

We search for exceptional collections whose $\Hom$- and $\Ext^1$-groups are mostly zero. In good cases, this implies that $B$ is formal, and $H^*B$ has no deformations. It then follows that $\eE$ is rigid, i.e.~constant in families.

\subsection{Quasiphantoms on the Kulikov surface}\label{sec!Kulikov-quasiphantom}
We study some properties of the exceptional collections on the Kulikov surface from Section \ref{sec!Kulikov}. For the purposes of the discussion, we fix the following exceptional collection
\[\EE\colon \Oh,\ L_1[2,2,0],\ L_2[2,1,2],\ L_3[0,0,1],\ L_4[1,1,1],\ L_5[2,2,1],\]
which can be found in the second row of Table \ref{tab!Kulikov-exceptional} in Section \ref{sec!Kulikov}.

Using Theorem \ref{thm!pushforward}, we may compute the $\Ext$-groups of the extended sequence \eqref{eqn!antican-extend}. We present the results in Table \ref{tab!ext-Kulikov} below. The $ij$th entry of the table is the following formal polynomial in $q$
 \[\sum_{k\in\ZZ}\dim\Ext^k(E_i,E_{i+j})q^k,\]
where $0\le i,j\le 5$, and the zigzag delineates those entries whose target $E_{i+j}$ is in the anticanonically extended part of \eqref{eqn!antican-extend}.
\begin{table}[ht]
\[
\renewcommand{\arraystretch}{1.3}
\begin{array}{c||cccccc}
&0&1&2&3&4&5\\
\hline\hline
0&1&2q^2&2q^2&2q^2&3q^2&\multicolumn{1}{c|}{3q^2}\\\cline{7-7}
1&1&0&0&2q+3q^2&q+2q^2&\multicolumn{1}{|c}{4q^2}\\\cline{6-6}
2&1&0&q^2&q^2&\multicolumn{1}{|c}{4q^2}&6q^2\\\cline{5-5}
3&1&q^2&q^2&\multicolumn{1}{|c}{4q^2}&6q^2&6q^2\\\cline{4-4}
4&1&0&\multicolumn{1}{|c}{3q^2}&5q^2&5q^2&5q^2\\\cline{3-3}
5&1&\multicolumn{1}{|c}{3q^2}&5q^2&5q^2&5q^2&6q^2\\\cline{2-2}
\end{array}\]
\caption{Ext-table of an exceptional collection on the Kulikov surface}\label{tab!ext-Kulikov}
\end{table}
\begin{lemma}\label{lem!Kulikov-ext1} The only nonzero $\Ext^1$-groups are $\Ext^1(E_1,E_4)$ which is $2$-dimensional, and $\Ext^1(E_1,E_5)$ which is $1$-dimensional.$\qed$
\end{lemma}
\begin{rmk}\rm The lemma shows that $\EE$ does not have $3$-block structure. A $3$-block structure means the exceptional collection can be split into three mutually orthogonal blocks (cf.~\cite{KN}). In fact, every exceptional collection in Table \ref{tab!Kulikov-exceptional}, and every exceptional collection in the Weyl group orbit (cf.~Section \ref{sec!group-actions}), has some non-zero $\Ext^1$-groups. This is in contrast with the exceptional collections on the Burniat surface exhibited in \cite{AO}, which are of the same numerical type, and have $3$-block structure.
\end{rmk}
\begin{prop} The $A_\infty$-algebra of $\EE$ is formal, and the product $m_2$ of any two elements with strictly positive degree is trivial.
\end{prop}
\begin{pf} The $A_\infty$-algebra $H^*B$ of $\EE$, is the direct sum of all $\Ext$-groups appearing above the zigzag in the table. By \cite[Lemma 2.1]{Sei} or \cite[Theorem 3.2.1.1]{Lef}, we may assume that $m_n(\dots,id_{E_i},\dots)=0$ for all $E_i$ and all $n>2$.

We show that every product $m_3$ must be zero for degree reasons. By Lemma \ref{lem!Kulikov-ext1}, there are only two nonzero arrows in degree $1$, and they can not be composed with one another, since they have the same source. Thus the product $m_3$ of any $3$ composable elements of $H^*B$ has degree at least $\deg m_3+1+2+2=4$, and is therefore identically zero, because the graded piece $H^4B$ is trivial. The same argument applies for all products $m_n$ with $n\ge 3$. Thus $H^*B$ is a formal $A_\infty$-algebra. In fact, we see from the table that any product $m_2$ of two elements of nonzero degree also vanishes for degree reasons.$\qed$
\end{pf}
Moreover, we calculate the Hochschild cohomology of $\Ay$ using heights.
\begin{prop} We have $\HH^0(\Ay)=\CC$, $\HH^1(\Ay)=0$, $\HH^2(\Ay)=\CC$, and $\HH^3(\Ay)$ contains a copy of $\CC^3$.
\end{prop}
\begin{pf} The pseudoheight of $\EE$ may also be computed from the table, where now we also need the portion below the zigzag. The minimal contribution to the pseudoheight is achieved by incorporating one of the nonzero $\Ext^1$-groups. For example,
\[e(E_1,E_4)+e(E_4,E_1\otimes\omega_X)-1+2=1+2-1+2=4,\] 
so $ph(\EE)=4$. In this case, by \cite{K}, the height and pseudoheight are equal. Hence $\HH^k(\Ay)=\HH^k(X)$ for $k\le2$, and $\HH^3(\Ay)\supset\HH^3(X)$. By the Hochschild--Kostant--Rosenberg isomorphism, the dimensions of $\HH^k(X)$ follow from the infinitesimal deformation theory of the Kulikov surface, which was studied in \cite{CC}: $H^1(T_X)=1$ and $H^2(T_X)=3$.$\qed$
\end{pf}
In summary, we have
\begin{thm} Every Kulikov surface $X$ has a semiorthogonal decomposition
\[D^b(X)=\left<\eE,\Ay\right>\]
where $\eE$ is generated by the exceptional collection $\EE$, and $\eE$ is rigid, i.e. $\eE$ does not vary with $X$. The semiorthogonal complement $\Ay$ is a quasiphantom category whose formal deformation space is isomorphic to that of $D^b(X)$, and therefore $X$ may be reconstructed from $\Ay$.
\end{thm}

\section{Secondary Burniat surfaces and effective divisors}\label{sec!burniat-secondary}
Burniat surfaces were discovered in \cite{Bur}, and an alternate construction is given in \cite{Inoue}. There are several cases $X_k$, with $K^2=k$ for $2\le k\le 6$. For details we refer to \cite{Peters}, \cite{BC0a}. Exceptional collections on primary Burniat surfaces $X_6$ with $K^2=6$ were first constructed and studied in \cite{AO}, where two $3$-block exceptional collections are exhibited. Burniat surfaces with $K^2=3,4,5,6$ can be constructed as abelian covers satisfying assumptions (A), and so we are able to enumerate exceptional collections on all these Burniat surfaces. We do not reproduce these computations here, but see \cite{computer-code}. Exceptional collections of line bundles of maximal length on the Burniat--Campedelli surface $X_2$ with $K^2=2$ remain elusive, because this surface does not satisfy assumption (A1).

In computing exceptional collections on fake del Pezzo surfaces, it becomes clear that a characterisation of effective line bundles is very useful. In this section we prove the following theorem for the secondary nodal Burniat surface.
\begin{thm}\label{thm!Burniat-effective} Let $X$ be a nodal secondary Burniat surface with $K^2=4$. Then the semigroup of effective divisors on $X$ is generated by the reduced pullbacks of irreducible components of the branch divisor, together with the pullbacks $E_4$, $E_5$ of two $(-1)$-curves on $Y$.
\end{thm}
With appropriate changes, the same proof works for the other surfaces satisfying assumptions (A). Indeed, Theorem \ref{thm!Kulikov-effective} above for the Kulikov surface is an easier case of this result. The additional complexity here arises from two sources: some of the exceptional curves on $Y$ are not branch divisors, and there is a $(-2)$-curve.

\subsection{Burniat surfaces revisited}\label{sec!Burniat}

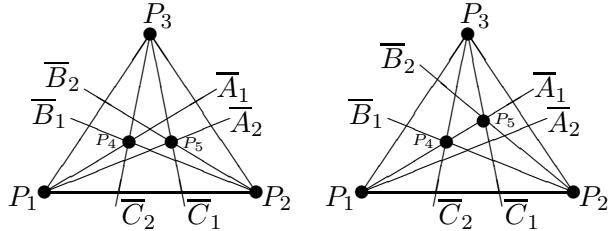
\begin{figure}[ht]
\begin{picture}(220,70)(-60,5)
\put(10,10){\line(1,0){80}}
\put(90,10){\line(-2,3){40}}
\put(10,10){\line(2,3){40}}
\put(10,10){\line(5,2){70}}
\put(90,10){\line(-5,3){65}}
\put(10,10){\line(5,3){65}}
\put(90,10){\line(-5,2){70}}
\put(50,70){\line(-1,-5){13}}
\put(50,70){\line(1,-5){13}}
\put(10,10){\circle*{5}}
\put(90,10){\circle*{5}}
\put(50,70){\circle*{5}}
\put(42,29){\circle*{5}}
\put(58,29){\circle*{5}}
\put(75,46){$\Abar_1$}
\put(80,33){$\Abar_2$}
\put(5,36){$\Bbar_1$}
\put(10,50){$\Bbar_2$}
\put(39,-2){$\Cbar_2$}
\put(64,-2){$\Cbar_1$}
\put(-4,5){$P_1$}
\put(92,5){$P_2$}
\put(47,74){$P_3$}
\put(29,28){\tiny{$P_4$}}
\put(62,27){\tiny{$P_5$}}
\put(130,10){\line(1,0){80}}
\put(210,10){\line(-2,3){40}}
\put(130,10){\line(2,3){40}}
\put(130,10){\line(5,2){70}}
\put(210,10){\line(-6,5){58}}
\put(130,10){\line(5,3){65}}
\put(210,10){\line(-5,2){70}}
\put(170,70){\line(-1,-5){13}}
\put(170,70){\line(1,-5){13}}
\put(130,10){\circle*{5}}
\put(210,10){\circle*{5}}
\put(170,70){\circle*{5}}
\put(162,29){\circle*{5}}
\put(176,37){\circle*{5}}
\put(195,46){$\Abar_1$}
\put(200,33){$\Abar_2$}
\put(125,36){$\Bbar_1$}
\put(137,57){$\Bbar_2$}
\put(159,-2){$\Cbar_2$}
\put(184,-2){$\Cbar_1$}
\put(116,5){$P_1$}
\put(212,5){$P_2$}
\put(167,74){$P_3$}
\put(149,28){\tiny{$P_4$}}
\put(180,36){\tiny{$P_5$}}
\end{picture}
\caption{The secondary Burniat configurations with $K^2=4$ (nodal configuration is on the right)}\label{fig!secondary-Burniat-4}
\end{figure}

We first describe the nodal secondary Burniat line configuration. Take the three coordinate points $P_1,P_2,P_3$ in $\PP^2$, and label the edges $\Abar_0=P_1P_2$, $\Bbar_0=P_2P_3$, $\Cbar_0=P_3P_1$. Then let $\Abar_1$, $\Abar_2$ (respectively $\Bbar_i$, $\Cbar_i$) be two lines passing through $P_1$ (resp.~$P_2$, $P_3$). We require that $\Abar_1$, $\Bbar_1$, $\Cbar_2$ are concurrent in $P_4$ (respectively $\Abar_1$, $\Bbar_2$, $\Cbar_1$ in $P_5$). This gives nine lines in total, four passing through each of $P_1$, $P_2$, $P_3$ and three passing through each of $P_4$, $P_5$. Moreover, $\Abar_1$ passes through three triple points. Blow up the five points $P_i$ to obtain a weak del Pezzo surface $Y$ of degree $4$. The strict transforms of these nine lines (for which we use the same labels) together with the three exceptional curves $\Ebar_i$ for $i=1,2,3$, are called the \emph{nodal secondary Burniat configuration} (see Figure \ref{fig!secondary-Burniat-4}).

The nodal secondary Burniat surface $X_4$ with $K_X^2=4$ is a $(\ZZ/2)^2$-cover of $Y$ branched in the configuration of Figure \ref{fig!secondary-Burniat-4}, and $X_4$ is a surface of general type with $p_g=0$, $K^2=4$ and $\Tors X=(\ZZ/2)^4$. The cover is not ramified in $\Ebar_4$, $\Ebar_5$. The weak del Pezzo surface $Y$ has a $(-2)$-curve, $\Abar_1$ and the canonical model of $X_4$ is a $(\ZZ/2)^2$-cover of a nodal quartic. Following the description of \cite{BC0b,BC0c}, the nodal secondary Burniat surfaces form an irreducible closed family, inside a $3$-dimensional irreducible connected component of the moduli space. This component is given by the union of the family of nodal secondary Burniat surfaces with the \emph{extended secondary Burniat surfaces}, which form an open subset. We do not directly consider extended Burniat surfaces here.

In Appendix \ref{app!Burniat-surface-data}, we show that the secondary Burniat surface satisfies assumptions (A). More precisely, we exhibit an explicit basis $e_0,\dots,e_5$ for $\Pic X/{\Tors X}$, in terms of reduced pullbacks of irreducible branch divisors. The appendix also lists coordinates for the reduced pullback of each irreducible component of the branch divisor, according to Definition \ref{def!coordinates}.

We define $\Eff$ to be the semigroup generated by the reduced pullbacks $A_0,\dots,C_2$, $E_1,E_2,E_3$ together with ordinary pullbacks $E_4$, $E_5$. There is a multigrading on $\Eff$ by multidegree in $\Pic X/{\Tors X}$, and we write $\Eff(d)$ for graded piece of multidegree $d$. Using the coordinates from Appendix \ref{app!Burniat-surface-data}, we define a homomorphism $t\colon\Eff\to\Tors X$, sending each generator of $\Eff$ to its associated torsion twist. Remember that $t(E_4)=t(E_5)=0$ because these are pulled back from $Y$.

\subsection{Proof of Theorem \ref{thm!Burniat-effective}.}

The stategy of proof is similar to that of Alexeev, \cite{Alexeev}, but for completeness, we outline the whole proof. The main differences are the $(-2)$-curve on $Y$ and the $(-1)$-curves which are not branch divisors. These introduce new complications which are not present in \cite{Alexeev}. We are able to resolve these issues because we can use the pushforward formula and Appendix \ref{app!Burniat-surface-data} to check effectivity in a systematic manner. 

Suppose $D$ is an effective divisor and $C$ is an effective curve class on $X$ for which $D\cdot C<0$. Then $C^2<0$ and $C$ is in the base locus of $D$, so we define $D'=D-aC$ where $a$ is the smallest positive integer for which $(D-aC)\cdot C\ge0$. In this way, we can reduce an effective divisor on $X$ to one which has positive intersection with all curve classes in $\Eff$. Such divisor classes form a rational polyhedral cone $\Pe$ in $N^1(X,\mathbb{R})$.

To describe the generators of $\Pe$ in the most geometric way, we first construct certain divisor classes on $X$ in terms of reduced pullbacks and the birational transformations the del Pezzo surface $Y$. Suppose we take a standard Cremona transformation of $\PP^2$ centred on any three non-collinear triple points $P_i$, $P_j$ and $P_k$. The numerical class of the hyperplane section of the image $\PP^2$ is $h_{ijk}=2e_0-e_i-e_j-e_k$ for any $\{i,j,k\}\ne\{1,4,5\}$, or $h_0=e_0$. There are also natural fibrations on $Y$ which arise from the pencil of hyperplanes passing through a fixed $P_k$ on some copy of $\PP^2$. The numerical classes of these fibrations are denoted $f_i=e_0-e_i$ or $f_{ijkl}=2e_0-e_i-e_j-e_k-e_l$ with $\{1,4,5\}\not\subset\{i,j,k,l\}$.

\begin{lemma}\label{lem!poly-gens} The polyhedron $\Pe$ is generated by the ten hyperplane classes $h_0$, $h_{ijk}$ and eight fibrations $f_i$, $f_{ijkl}$ defined above, together with the four additional classes
\begin{align*}
g_1&=3e_0-e_1-2e_2-e_3-e_4-e_5,&g_2&=3e_0-e_1-e_2-2e_3-e_4-e_5,\\
g_3&=3e_0-e_1-2e_2-e_4-e_5,& g_4&=3e_0-e_1-2e_3-e_4-e_5.
\end{align*}
\end{lemma}
\begin{pf} Any generator $D$ of $\Eff$ determines a linear function $\cdot D$, which in turn defines a collection of hyperplanes supporting the polyhedron $\Pe$. We use the computer \cite{computer-code} to calculate the integral generators of the cone.

We examine the additional generators. The class $g_1$ is the hyperplane section of the copy of $\PP^2$ obtained by contracting $\Abar_0$, $\Abar_1$, $\Bbar_0$, $\Bbar_1$ and $\Bbar_2$ on $Y$, and $g_3$ is the hyperplane section of the quadric cone given by contracting $\Abar_0$, $\Abar_1$, $\Ebar_3$, $\Bbar_1$ and $\Bbar_2$. There are similar descriptions of $g_2$ and $g_4$.
$\qed$
\end{pf}

\begin{lemma} Suppose $D$ is an effective divisor on $X$ with $K_X\cdot D\le4$. Then $D$ is in $\Eff$.
\end{lemma}
\begin{pf} We may assume that $D$ is in $\Pe$. This is a finite (and small) number of classes to check, and we do this directly using the computer implementation \cite{computer-code} of our pushforward formula Theorem \ref{thm!pushforward}.$\qed$
\end{pf}

\begin{prop}\label{prop!effective-chi-pos} Suppose $D$ is an effective divisor on $X$ with $K_X\cdot D>4$ and $\chi(D)>0$. Then $D$ is in $\Eff$.
\end{prop}
\begin{pf} Since $K_X\cdot D>4$ we have that $(K_X-D)\cdot K_X<0$ and so $K_X-D$ can not be effective. By Serre duality, $h^2(D)=h^0(K_X-D)=0$, hence $D$ is effective.

Choose $\bar{D}$ in $\Pic Y$ such that the numerical class of $K_Y+\bar{D}$ in $\Pic Y$ is the same as that of $D-K_X$ in $\Pic X/{\Tors X}$ by assumption (A). Then
\[\chi(K_Y+\bar{D})=1+\frac12(K_Y+\bar{D})\bar{D}=1+\frac12(D-K_X)D=\chi(D)>0.\]
Moreover, $h^2(K_Y+\bar{D})=h^0(-\bar{D})=0$ because ${-K_Y}\cdot{-\bar{D}}<0$, so by the same argument as above, we see that $K_Y+\bar{D}$ is effective on $Y$.

Now, any effective divisor on $Y$ is a positive linear combination of branch divisors $\Abar_0,\dots,\Cbar_2,\Ebar_1,\Ebar_2,\Ebar_3$ and exceptional curves $\Ebar_4$ and $\Ebar_5$. So taking the reduced pullback, we get the following expression for the numerical class of $D$ in $\Pic X/{\Tors X}$:
\[D = K_X+(\text{combination of }A_0,\dots,E_3)+\tfrac12(\text{combination of }E_4,E_5).\]
The coefficient of $\frac12$ appears because $\Ebar_4$ and $\Ebar_5$ are not branch divisors.
It remains to show that $D+\tau$ is in $\Eff$ for any $\tau$ such that $D+\tau$ is effective. This is implied by the following lemma:
\end{pf}
\begin{lemma}\label{lem!effective-chi-gt-0}
\begin{enumerate}
\item[(1)] Let $L$ be any of the following line bundles on $X$:
\begin{gather*}
\Oh_X(K_X+\gam)\otimes\tau,\ \Oh_X(K_X+\tfrac12E_4)\otimes\tau,\ \Oh_X(K_X+\tfrac12E_5)\otimes\tau,\text{ or }\\
\Oh_X(K_X+\tfrac12(E_4+E_5))\otimes\tau
\end{gather*}
where $\gam$ is any generator of $\Eff$ and $\tau$ is any element of $\Tors X$. Then $L$ is effective and in $\Eff$ unless $L=\Oh_X(K_X+A_1)$.
\item[(2)] The line bundles $L=\Oh_X(K_X+kA_1)$ are not effective for any $k>0$.
\end{enumerate}
\end{lemma}
\begin{pf}
\paragraph{(1)} Suppose $L=\Oh_X(K_X+A_0)\otimes\tau$, and take the graded piece of $\Eff$ with multidegree $d=(4,-2,-2,-1,-1,-1)$. We use the computer \cite{computer-code} to check that the image of $\Eff(d)$ under $t$ is all of $\Tors X$. This proves that $L$ is effective and in $\Eff$ for any $\tau$. The same computation works for all multidegrees listed in the statement, except when $L=\Oh_X(K_X+A_1)\otimes\tau$, for which we refer to the proof of part (2).
\paragraph{(2)} When $L=\Oh_X(K_X+A_1)\otimes\tau$, the same computation as above shows that the image of $\Eff(4,-2,-1,-1,-2,-2)$ under $t$ is $\Tors X-\{[1,0,0,0]\}$. Thus $\Oh_X(K_X+A_1)$ is not in $\Eff$. Indeed, the pushforward is
\[\fie_*L=\Oh_Y(e_2-e_1)\oplus
\Oh_Y(e_3-e_2)\oplus\Oh_Y(-2e_0+e_2+e_3)\oplus\Oh_Y(e_0-e_3-e_4-e_5),\]
which is not effective. Moreover, by the projection formula, we have
\[\fie_*L(2mA_2)=\fie_*L\otimes\Oh_Y(m\Abar_2)=\fie_*L\otimes\Oh_Y(m(e_0-e_1-e_4-e_5)),\]
which is not effective for any $m$, and so $\Oh_X(K_X+kA_2)$ is not effective for any odd $k=1+2m$. For even $k$, the proof is similar, starting from $\fie_*\Oh_X(K_X)$.
$\qed$
\end{pf}
\begin{rmk}\rm Since $A_2$ is a $(-2)$-curve, we have $K_X\cdot(K_X+kA_1)=K_X^2=4$ for all $k$. Thus we do not need part (2) of the above lemma, because $K_X+kA_1$ does not satisfy the assumptions of Proposition \ref{prop!effective-chi-pos}.
\end{rmk}
Finally, we take care of the cases with $\chi(D)\le0$.
\begin{lemma}\label{lem!chi-zero} Suppose $D$ is an effective divisor on $X$ with numerical class in $\Pe$ and $\chi(D)\le0$. Then $D$ is in one of the following classes:
\begin{enumerate}
\item[(1)] $h$ or $2h$ for any hyperplane generator $h$;
\item[(2)] $g$ or $2g_1$ or $2g_3$, where $g$ refers to any of the additional generators described in Lemma \ref{lem!poly-gens};
\item[(3)] $nf$, $nf+f'$, $nf+h$, $nf+g$ for any $n\ge1$ where $f$ is a fibration and $f'$ is another fibration with intersection $f\cdot f'=1$, $h$ is a hyperplane generator with $f\cdot h=1$, $g$ is an additional generator with $f\cdot g=1$.
\end{enumerate}
\end{lemma}
\begin{pf} This is a systematic induction. We note that each generator $\gam$ of $\Pe$ has $\chi(\gam)=0$. Moreover, if $D=D_1+D_2$ then $\chi(D)=\chi(D_1)+\chi(D_2)+D_1\cdot D_2-1$. So for example, starting from $f_1$, we choose another fibration generator $f'$. Either $f_1\cdot f'=0$, in which case $f_1=f'$ and $\chi(2f_1)=-1$, or $f_1\cdot f'=1$, so that $\chi(f_1+f')=0$. Now adding a further generator $\gam$ to $f_1+f'$ yields $\chi(f_1+f'+\gam)>0$ by simple consideration of the intersection numbers, unless $\gam$ is one of $f_1$ or $f'$. We continue in this way, to produce the list of possibilities.$\qed$
\end{pf}
\begin{lemma} Suppose $L$ is an effective line bundle with numerical class one of the exceptional cases from Lemma \ref{lem!chi-zero}. Then $L$ is in $\Eff$.
\end{lemma}
\begin{pf} We give a proof for $nf_1$. The other possibilities listed in Lemma \ref{lem!chi-zero}(3) work in the same way, and cases (1) and (2) can be checked by a direct computation \cite{computer-code}. As in the proof of Lemma \ref{lem!effective-chi-gt-0}, we split into even and odd cases and make use of the projection formula. 

Let $L=\Oh_X(2f_1)\otimes\tau$ for some $\tau$ in the image of $t(\Eff(2f_1))$, so that in particular, $L$ is effective. Since $C_0+E_3$ is a section of $\Oh_X(f_1)$, it follows that $\Oh_X(nf_1)\otimes\tau$ is effective and in $\Eff$ for any $n\ge2$.

Now suppose $\tau$ is any torsion element in $\Tors X - t(\Eff(2f_1))$, so that $L$ is not in $\Eff$. For example, $\tau=[0,0,0,1]$. Then we compute
\begin{align*}
\fie_*L=\Oh_Y( 0, -1,  0,  0,  0,  1)&\oplus \Oh_Y(-2,  1,  1,  1,1,1)\\
&\oplus\Oh_Y(-1, -1,  1,  1,  1,  0)\oplus\Oh_Y(-2,  0,  1,  1,  1,  1),
\end{align*}
which is clearly not effective. Moreover, by the projection formula, we see that $\fie_*L\otimes\Oh_X(2mf_1)=\fie_*L\otimes\Oh_Y(mf_1)$ is not effective for any $m$ either, for degree reasons. This completes the proof for any even multiple of $f_1$. A similar computation proves the odd case, starting from $3f_1$.$\qed$
\end{pf}
%\vfill
%\bigskip\bigskip\bigskip\bigskip\bigskip\bigskip
%%%%%%%%%%%%%%%%%%%%%%%%%%%%%%%%%%%%%%%%%%%%%%%%%%%%%%%%%%%%%%%%%%%%%%
\appendix

\section{Appendix: Acyclic bundles on the Kulikov surface}\label{app!Kulikov}
For reference, here are the acyclic line bundles on the Kulikov surface used in Section \ref{sec!Kulikov}.
\[\renewcommand{\arraystretch}{1.3}
\begin{array}{c||l}
L&\Ay(L)\\
\hline\hline
L_1^{-1}&[ 0, 0, 0 ], [ 0, 1, 0 ], [ 1, 1, 0 ], [ 2, 1, 0 ], [ 2, 2, 0 ], [ 1, 0, 1 ], [ 2, 0, 1 ], [ 0, 1, 1 ],\\
& [ 1, 1, 1 ], [ 0, 2, 1 ], [ 2, 2, 1 ], [ 0, 1, 2 ], [ 1, 1, 2 ], [ 0, 2, 2 ]\\
\hline
L_2^{-1}&[ 0, 1, 0 ], [ 1, 1, 0 ], [ 2, 2, 0 ], [ 2, 0, 1 ], [ 0, 1, 1 ], [ 1, 1, 1 ], [ 2, 1, 1 ], [ 1, 2, 1 ],\\
& [ 2, 2, 1 ], [ 0, 0, 2 ], [ 1, 0, 2 ], [ 0, 1, 2 ], [ 1, 1, 2 ], [ 1, 2, 2 ]\\
\hline
L_3^{-1}&[ 0, 1, 0 ], [ 1, 1, 0 ], [ 1, 0, 1 ], [ 0, 1, 1 ], [ 1, 1, 1 ], [ 0, 2, 1 ], [ 1, 2, 1 ], [ 0, 0, 2 ],\\
& [ 2, 0, 2 ], [ 0, 1, 2 ], [ 1, 1, 2 ], [ 2, 1, 2 ], [ 0, 2, 2 ], [ 1, 2, 2 ]\\
\hline
L_4^{-1}&[ 0, 0, 0 ], [ 0, 1, 0 ], [ 2, 1, 0 ], [ 0, 2, 0 ], [ 2, 2, 0 ], [ 1, 0, 1 ], [ 2, 0, 1 ], [ 0, 1, 1 ],\\
& [ 1, 1, 1 ], [ 2, 1, 1 ], [ 0, 2, 1 ], [ 2, 2, 1 ], [ 1, 1, 2 ], [ 0, 2, 2 ], [ 2, 2, 2 ]\\
\hline
L_5^{-1}&[ 0, 1, 0 ], [ 1, 1, 0 ], [ 2, 2, 0 ], [ 1, 0, 1 ], [ 2, 0, 1 ], [ 0, 1, 1 ], [ 1, 1, 1 ], [ 0, 2, 1 ],\\
& [ 1, 2, 1 ], [ 2, 2, 1 ], [ 0, 0, 2 ], [ 0, 1, 2 ], [ 1, 1, 2 ], [ 0, 2, 2 ], [ 1, 2, 2 ]\\
\hline\hline
L_2^{-1}\otimes L_1&[ 1, 0, 0 ], [ 2, 0, 0 ], [ 2, 1, 0 ], [ 0, 1, 1 ], [ 0, 1, 2 ], [ 2, 1, 2 ], [ 0, 2, 2 ]\\
\hline
L_3^{-1}\otimes L_1&[ 0, 0, 0 ], [ 1, 0, 0 ], [ 2, 0, 0 ], [ 1, 1, 0 ], [ 2, 1, 0 ], [ 2, 2, 0 ], [ 1, 1, 2 ], [ 2, 1, 2 ],\\
& [ 2, 2, 2 ]\\
\hline
L_4^{-1}\otimes L_1&[ 0, 1, 0 ], [ 1, 1, 0 ], [ 0, 1, 1 ], [ 1, 1, 1 ], [ 1, 2, 1 ], [ 0, 0, 2 ], [ 1, 0, 2 ], [ 2, 0, 2 ],\\
& [ 0, 1, 2 ], [ 1, 1, 2 ], [ 1, 2, 2 ]\\
\hline
L_5^{-1}\otimes L_1&[ 1, 0, 0 ], [ 2, 0, 0 ], [ 1, 1, 0 ], [ 2, 1, 0 ], [ 2, 2, 0 ], [ 0, 1, 1 ], [ 0, 0, 2 ], [ 0, 1, 2 ],\\
& [ 1, 1, 2 ], [ 2, 1, 2 ], [ 0, 2, 2 ], [ 2, 2, 2 ]\\
\hline
L_3^{-1}\otimes L_2&[ 1, 0, 1 ], [ 1, 1, 1 ], [ 2, 1, 1 ], [ 2, 0, 2 ], [ 1, 1, 2 ], [ 2, 1, 2 ], [ 1, 2, 2 ]\\
\hline
L_4^{-1}\otimes L_2&[ 0, 0, 0 ], [ 0, 1, 0 ], [ 1, 1, 0 ], [ 1, 0, 1 ], [ 0, 1, 1 ], [ 1, 1, 1 ], [ 0, 2, 1 ], [ 2, 0, 2 ],\\
& [ 0, 1, 2 ], [ 1, 1, 2 ], [ 0, 2, 2 ]\\
\hline
L_5^{-1}\otimes L_2&[ 0, 1, 0 ], [ 1, 0, 1 ], [ 0, 1, 1 ], [ 1, 1, 1 ], [ 2, 1, 1 ], [ 0, 2, 1 ], [ 0, 0, 2 ], [ 2, 0, 2 ],\\
&[ 1, 1, 2 ], [ 2, 1, 2 ], [ 0, 2, 2 ], [ 1, 2, 2 ]\\
\hline
L_4^{-1}\otimes L_3&[ 0, 0, 0 ], [ 0, 1, 0 ], [ 1, 1, 0 ], [ 2, 2, 0 ], [ 2, 0, 1 ], [ 0, 1, 1 ], [ 1, 1, 1 ], [ 2, 2, 1 ],\\
& [ 1, 0, 2 ], [ 0, 1, 2 ], [ 1, 1, 2 ]\\
\hline
L_5^{-1}\otimes L_3&[ 0, 1, 0 ], [ 1, 1, 0 ], [ 2, 1, 0 ], [ 2, 2, 0 ], [ 2, 0, 1 ], [ 1, 1, 1 ], [ 2, 1, 1 ], [ 1, 2, 1 ],\\
& [ 0, 0, 2 ], [ 1, 0, 2 ], [ 0, 1, 2 ], [ 1, 2, 2 ]\\
\hline
L_5^{-1}\otimes L_4&[ 1, 0, 0 ], [ 2, 0, 0 ], [ 1, 1, 0 ], [ 2, 2, 0 ], [ 0, 0, 2 ], [ 0, 1, 2 ], [ 2, 1, 2 ],[ 2, 2, 2 ]
\end{array}\]

\pagebreak

\section{Appendix: Nodal Secondary Burniat surface with $K^2=4$}\label{app!Burniat-surface-data}
The maps $\Psi_4,\Psi_4^n\colon H_1(Y-\Del,\ZZ)\to(\ZZ/2)^6$ determining respectively the non-nodal and nodal Burniat surfaces, differ from one another slightly. We tabulate them below.
\[\renewcommand{\arraystretch}{1.3}\begin{array}{c||c|c|c|c|c|c|c|c|c}
\Gam&\Abar_0&\Abar_1&\Abar_2&\Bbar_0&\Bbar_1&\Bbar_2&\Cbar_0&\Cbar_1&\Cbar_2\\
\hline\hline
\Psi_4(\Gam)-\Fie(\Gam)&0&g_3&g_4&0&g_5&g_6&0&g_4+g_6&g_3+g_5\\
\Psi_4^n(\Gam)-\Fie(\Gam)&0&g_3&g_4&0&g_5&g_6&g_3+g_4&g_3+g_6&g_3+g_5
\end{array}\]
The restriction imposed by $P_5$ is $\Psi_4(\Abar_2+\Bbar_2+\Cbar_1)=0$ in the non-nodal case, and $\Psi^n_4(\Abar_1+\Bbar_2+\Cbar_1)=0$ in the nodal case. Either way, $g_7$ is eliminated, so the torsion group is $(\ZZ/2)^4$, generated by $g_3^*,\dots,g_6^*$.

We extend the basis chosen for the free part of $\Pic(X_5)$. The basis is the same for non-nodal and nodal surfaces
\begin{gather*}e_0=C_0+E_1+E_3,\ e_1=E_1,\ e_2=E_2,\ e_3=E_3,\\
e_4=C_0-C_2+E_1,\ e_5=B_0-B_2+E_3.
\end{gather*}
Coordinates for non-nodal surface:
\[\renewcommand{\arraystretch}{1.1}\begin{array}{c|rrrrrr|c}
&\multicolumn{6}{c|}{\text{Multidegree}}&\text{Torsion}\\
\hline\hline
\Oh_X(A_0)&1&-1&-1&0&0&0&[1,1,0,0]\\
\Oh_X(A_1)&1&-1&0&0&-1&0&[1,0,0,0]\\
\Oh_X(A_2)&1&-1&0&0&0&-1&[0,1,1,0]\\
\Oh_X(B_0)&1&0&-1&-1&0&0&[0,0,1,1]\\
\Oh_X(B_1)&1&0&-1&0&-1&0&[0,0,1,0]\\
\Oh_X(B_2)&1&0&-1&0&0&-1&[0,0,1,1]\\
\Oh_X(C_0)&1&-1&0&-1&0&0&0\\
\Oh_X(C_1)&1&0&0&-1&0&-1&[0,0,1,0]\\
\Oh_X(C_2)&1&0&0&-1&-1&0&0
\end{array}\]
Coordinates for nodal surface are the same (with same multidegrees) except for the following:
\[\renewcommand{\arraystretch}{1.1}\begin{array}{c|rrrrrr|c}
&\multicolumn{6}{c|}{\text{Multidegree}}&\text{Torsion}\\
\hline\hline
%\Oh_X(D_1)&1&-1&-1&0&0&0&[1,1,0,0]\\
\Oh_X(A_1)&1&-1&0&0&-1&-1&[1,0,1,0]\\
\Oh_X(A_2)&1&-1&0&0&0&0&[0,1,0,0]\\
%\Oh_X(D_2)&1&0&-1&-1&0&0&[0,0,1,1]\\
%\Oh_X(D_2^1)&1&0&-1&0&0&-1&[0,0,1,1]\\
%\Oh_X(D_2^2)&1&0&-1&0&-1&0&[0,0,0,1]\\
%\Oh_X(D_3)&1&-1&0&-1&0&0&0\\
%\Oh_X(D_3^1)&1&0&0&-1&0&-1&[0,0,0,1]\\
%\Oh_X(D_3^2)&1&0&0&-1&-1&0&0
\end{array}\]
In both cases, $\Oh_X(K_X)=\Oh(3,-1,-1,-1,-1,-1)[0,0,1,0]$.
\vfill
%%%%%%%%%%%%%%%%%%%%%%%%%%%%%%%%%%%%%%%%%%%%%%%%%%%%%%%%%%%%%%%%

%%%%%%%%%%%%%%%%%%%%%%%%%%%%%%%%%%%%%%%%%%%%%%%%%%%%%%%%%%%%%%%%

%%%%%%%%%%%%%%%%%%%%%%%%%%%%%%%%%%%%%%%%%%%%%%%%%%%%%%%%%%%%%%%%%%%
\pagebreak

\end{document}